\documentclass[11pt]{amsart}

\usepackage{amsmath,amssymb,graphicx,xspace}

\usepackage{amsthm, ifthen, amsfonts, epsfig, graphics, latexsym, psfrag, verbatim, textcomp, tikz, slashbox}

\renewcommand{\labelenumi}{\arabic{enumi}.}
\renewcommand{\labelenumii}{\alph{enumii}.}

\newtheorem{thm}{Theorem}[section]
\newtheorem{lem}[thm]{Lemma}
\newtheorem{cor}[thm]{Corollary}

\newtheorem{prop}[thm]{Proposition}
\theoremstyle{definition}
\newtheorem{exmp}[thm]{Example}
\newtheorem{defn}[thm]{Definition}

\newcommand{\bbb}[1]{\ensuremath{\mathbb{#1}}}
\newcommand{\D}{\bbb{D}}
\newcommand{\euc}{\bbb{E}}
\newcommand{\R}{\bbb{R}}
\newcommand{\sph}{\bbb{S}}
\newcommand{\Z}{\bbb{Z}}
\newcommand{\cayley}{\textsc{Cay}}
\newcommand{\cat}{\textsc{CAT}}
\newcommand{\E}{\ensuremath{\mathcal{E}}}
\newcommand{\G}{\ensuremath{\mathcal{G}}}
\newcommand{\pe}{\textsf{PE}\xspace}
\newcommand{\ps}{\textsf{PS}\xspace}
\newcommand{\size}[1]{\ensuremath{\vert #1 \vert}}

\begin{document}

\title{Full Subcomplexes of  $\cat(0)$ Simplicial Complexes}

\author[R.~Levitt]{Rena M.H. Levitt}
\address{Department of Mathematics\\
        Pomona College\\
         610 North College Avenue\\
Claremont, CA 91711}
\email{rena.levitt@pomona.edu}

\date{\today}

\begin{abstract} 
In this article, I prove that full subcomplexes of $\cat(0)$ simplicial $3$-complexes inherit the non positive curvature condition, and describe a family of counterexamples that prove   this result can not be extended to higher dimensions.
\end{abstract}

\maketitle

\section{Introduction}

Since Gromov coined the term $\cat(0)$ is his  seminal 1987 paper  \cite{Gr87}, the $\cat(0)$ and nonpositively curvature conditions have been exceptionally useful in studying both the metric spaces that satisfy these conditions and the groups acting on them. For instance if a metric space $X$ is non positively curved, then its universal cover $\tilde{X}$ is contractible and the word problem for its fundamental group $\pi_1(X)$ can be solved in quadratic time.  Unfortunately showing that a space is $\cat(0)$ can in practice be quite difficult. For example, while Elder and McCammond proved the existence of an algorithm to determine if a finite metric polyhedral complex  is nonpositively curved \cite{ElMc04}, the algorithm they developed is theoretical. Few practical algorithms exist to determine if a given metric simplicial complex is nonpositively curved. 

In this article I focus on answering the following question; if $L$ is full subcomplex of a $\cat(0)$ simplicial complex $K$, does $L$ inherit nonpositive curvature condition from $K$? A subcomplex $L$ is \emph{full} in $K$ if $L$ contains each simplex of $K$ whose vertex set is contained in $L$. The following theorem is the main result of this article.

\newtheorem*{main}{Main Theorem \ref{thm:main}}

\begin{main}[Full $\Rightarrow$ NPC]\label{thm:main} Let $L$ be a full subcomplex of a regular $\cat(0)$ simplicial $3$-complex $K$. Then $L$ is nonpositively curved.
\end{main}

\noindent The proof of the main theorem follows from a careful analysis the structure of links in in full subcomplexes. It relies on two previously known results; Katherine Crowley's proof of the existence of $\cat(0)$ spanning disk  \cite{Cr08}, and Elder and McCammond's list of forbidden link configurations in $\cat(0)$ simplicial $3$-complexes  \cite{ElMc02}. 

Unfortunately Main Theorem \ref{thm:main} does not extend to higher dimensions. There exists a family of counterexamples to the main theorem in each  dimension greater than or equal to four. Excluding these counterexamples as subcomplexes gives complexes that fall under Januszkiewicz and \'{S}wi\c{a}tkowski's theory of simplicial non positive curvature \cite{JaSw06}.

The article is structured as follows. Section \ref{sec:prelim} is a brief review of curvature in simplicial complexes. The second section discusses combinatorial paths and  filling disks in simplicial $3$-complexes. In particular, Crowley's existence theorem is stated in Section \ref{sec:disks}, along with some useful corollaries of her result. Section \ref{sec:links} is devoted to the structure of metric links when restricted full subcomplexes. The proof of the main theorem is the subject of Section \ref{sec:proof}. Finally in Section \ref{sec:highdim} is a description of a family of counterexamples to the main theorem in higher dimensions.  


\section{Curvature in Simplicial Complexes}\label{sec:prelim}


This section is a brief review of piecewise euclidean and spherical simplicial complexes, and non positive curvature. Recall that a metric space is \emph{geodesic} if every pair of points is connected by a length minimizing path. 

\begin{defn}[Piecewise Euclidean and Spherical Complexes] A \emph{Euclidean polytope} is the convex hull of a finite set of points in euclidean space $\euc^{n}$. Similarly, a \emph{Spherical polytope} is the convex hull of a finite set of points contained in an open hemisphere of $\sph^{n}$. A \emph{piecewise Euclidean complex}, or \emph{PE-complex} is a cell complex built out of Euclidean polytopes glued together along faces by isometries.  A \emph{piecewise Spherical complex (PS-complex)} is a cell complex built out of spherical polytopes. A theorem of Martin Bridson's implies that piecewise Euclidean and Spherical complexes with finitely many cell isometry types are geodesic metric spaces \cite{BrHa99}. The dimension of a complex $K$ is the maximum of the dimensions of its cells if the maximum exists. If not, $K$ is infinite dimensional. The $n$-skeleton of a complex $K$, denoted $K^{(n)}$ is the union of $m$-cells of $K$ for $m \leq n$. A \emph{subcomplex} $L$ of a complex $K$ is a subset of $K$ that is also a complex. \end{defn} 

\begin{defn}[Simplicial Complex] An \emph{$n$-dimensional simplex} or \emph{$n$-simplex}  $\sigma$ is the convex hull of $n+1$ points in general linear position in $\euc^n$. A face of a simplex is the convex hull of a subset of the points defining $\sigma$. A $0$-simplex is a vertex, a $1$-simplex an edge, and a $2$-simplex a face. It can be useful to identify a simplex with its vertex set. If $\sigma^{(0)} = \{v_0, v_1, ..., v_n\}$, then $\{v_0, v_1, ..., v_n\}s$ \emph{span} $\sigma$. A \emph{simplicial complex} is a piecewise Euclidean complex with each $n$-polytope isometric to an $n$-simplex. A simplicial complex is \emph{regular} if each edge has unit length. All simplical complexes in this paper are taken to be regular unless otherwise specified. A \emph{spherical simplicial complex} is a piecewise spherical complex with each cell isometric to a simplex in an open hemisphere of $\sph^{n}$. 
\end{defn}

 \begin{defn}[Flag and Full] A simplicial complex $K$ is \emph{flag} if every set of vertices pairwise connected by edges span a simplex of $K$. A subcomplex $L$ of a complex $K$ of is \emph{full} in $K$ if  $\sigma^{(0)} \subseteq L$ implies $\sigma \subseteq L$ for all $\sigma \subseteq K$. It immediately follows that Full subcomplexes of flag complexes are flag. \end{defn}

Intuitively, a $\cat(0)$ space is a geodesic metric with geodesic triangles ``thinner'' than their Euclidean counter parts.  For the purposes of this paper $\cat(0)$ will be defined in terms of Gromov's link condition. This requires the following definitions. 

\begin{defn}[Metric Link] Let $\sigma$ be a $k$-face of an  $n$-simplex $\tau$. The \emph{metric link} of $\sigma$ in $\tau$ is  the set of unit tangent vectors orthogonal to $\sigma$ and pointing into $\tau$. This defines a spherical $(n-k-1)$-simplex. Let $\sigma$ be a cell of a simplicial complex $K$. The \emph{metric link} of $\sigma$ in $K$, denoted $lk_{K}(\sigma)$, is the spherical simplicial complex whose cells are the links of $\sigma$ in each $\tau \supseteq \sigma$ in $K$.  
\end{defn}

\begin{lem}\label{lem:linksfull} Let $L$ be a full subcomplex of a simplicial complex $K$. Then $lk_L(v)$ is full in $lk_K(v)$ for each vertex $v \subseteq L$. \end{lem}

\begin{proof} Suppose that $\sigma^{(0)} = \{v_0,v_1,...v_n\}\subseteq lk_L(v)$ for $\sigma \subseteq lk_K(v)$. Each vertex $v_i$ corresponds to an edge $e_i$ in $L$. Let $w_i$ be the vertex opposite $v$ on $e_i$ for $i = 0...n$. The existence of $\sigma \subseteq lk_K(v)$ implies that $\{v,w_1,w_2,...,w_n\}$ span a simplex $\tau$ in $K$. By fullness $\tau$ is contained in $L$. Thus $\sigma \subseteq lk_L(v)$.
\end{proof}

\begin{defn}[Nonpositively Curved, $\cat(0)$] Let $X$ be a geodesic metric space. A \emph{locally geodesic loop} $\rho$ is an embedding of a metric circle into $X$ satisfying the following property; at each point $x$ on the image of $\rho$ the angle between the incoming and outgoing tangent vectors of $\rho$ in $lk_{X}(x)$ is at least $\pi$. A piecewise Euclidean complex with finitely many isometry types of cells is \emph{nonpositively  curved} if the link of each cell contains no locally geodesic loops of length less than $2\pi$. If in addition $X$ is connected and simply connected, then $X$ is $\cat(0)$. \end{defn}

The intrinsic metric on a $\cat(0)$ space is convex, and geodesics are unique. 


\section{Combinatorial Paths and Filling Disks}\label{sec:disks}


This section is a discussion combinatorial paths and disks in $\cat(0)$ simplicial complexes, including Katherine Crowley's result for spanning disks in $3$-complexes. The end of the section includes some useful consequences of her result.

\begin{defn}[Combinatorial Path, and Combinatorial Distance] A \emph{combinatorial path} is an alternating sequence of vertices and edges  $\alpha = [v_0,e_1,v_1, \\ e_2,..., v_{n-1},e_n, v_{n}]$ such that $e_i$ is spanned $\{v_{i-1},v_{i}\}$ for $1 \leq i <k$. The path is a \emph{loop} if $v_0 = v_n$. A combinatorial path is \emph{tight} if it does not cross the same edge twice. The length of a combinatorial path $\ell(\alpha)$ is the number of edges it crosses. The \emph{combinatorial distance} between two vertices $v$ and $w$,  denoted $d_c(v,w)$, is the minimum of the set of lengths of combinatorial path from $v$ to $w$. This defines a metric on $K^{(0)}$. A combinatorial path $\gamma$  from $v$ to $w$ is a \emph{combinatorial geodesic} if $\ell(\gamma) = d_c(v,w)$. Combinatorial geodesics are not unique.\end{defn}

\begin{defn}[Disk Diagram] A \emph{disk diagram} is a contractible $2$-complex that embeds in $\euc^2$. This embedding is often implicit. A disk diagram is \emph{nonsingular} if it is homeomorphic to the closed unit disk. Otherwise the disk contains a \emph{cut point} whose removal disconnects the diagram. Disks with cut points are called \emph{singular}. Singular disks may broken up into nonsingular subdisks. The boundary of a disk diagram is a combinatorial loop read clockwise around the outside of the disk. This may be ambiguous for a singular disk, where the boundary is determined by giving an explicit embedding in the plane. Suppose $D$ is a simplicial disk, i.e. $D$ has triangular faces. The \emph{combinatorial area} of $D$ is the number of faces contained in $D$.\end{defn}

\begin{defn}[Spanning Disk] Let $D$ be a simplicial disk diagram, $K$ be a simplicial complex, and $f: D \longrightarrow K$ a cellular map. Then $\alpha= f(\partial D)$  is a combinatorial loop in  $K$. In this case, $D$ \emph{spans} $\alpha$ in $K$ and $\alpha$ bounds $D$. \end{defn}

\begin{defn}[Vertex Degree] Let $D$ be a simplicial disk diagram. The \emph{degree} of a vertex $v$ is the number of edges sharing $v$ as a vertex. For vertices on the interior of $D$, this is equivalent to the number of faces with $v$ as a vertex. For boundary vertices, this is one more than the number of faces with $v$ as a vertex. By Gromov's link condition a disk is $\cat(0)$ if and only if each interior vertex is contained in at least six triangles. \end{defn}

Many of the results in this paper involve analyzing the structure of spanning disks. The Combinatorial Gau\ss-Bonnet Theorem will be used in many of these arguments. This is a classical result, a proof can be found in \cite{Cr08} 

\begin{thm}[Combinatorial Gau\ss-Bonnet] Let $D$ be a triangulated disk. Then
\begin{equation*}
    \sum_{v \in \partial D} (4-deg(v)) + \sum_{v \in intD} (6-deg(v)) = 6.
\end{equation*}
\end{thm}

By contractibility, each combinatorial loop bounds a spanning disk. The following theorem of Katherine Crowley's  determines the structure of these disks in $\cat(0)$ simplicial $3$-complexes \cite{Cr08}.

\begin{thm}[Crowley, Spanning Disks] \label{thm:spanningdisks} Let $K$ be a $\cat(0)$ simplicial $3$-complex, $\alpha$ be a combinatorial loop in $K^{(1)}$. Then there exists a $\cat(0)$ disk $D$  contained in $K$ of minimal combinatorial area such that $\partial D = \alpha$.
\end{thm}

\begin{cor} Let $D$ be a minimal minimal spanning disk. Then  \begin{equation*} \sum_{v \in \partial D} (4-deg(v)) \geq 6 \end{equation*}. \end{cor}

\begin{proof}Theorem \ref{thm:spanningdisks} implies that the degrees of each interior vertex is at least six. Thus $\sum_{v \in intD} (6-deg(v)) <0$. \end{proof}

\begin{defn}[Empty $n$-gon]  A \emph{combinatorial $n$-gon} or \emph{$n$-gon} is a tight combinatorial loop of length $n$.  A $n$-gon $\alpha$ is \emph{empty} if the minimal disk spanning $\alpha$ has an interior vertex. For example, an empty triangle is a loop of length three not spanned by a face, an empty square is a loop of length four not spanned by two faces sharing an edge, and an empty pentagon is a loop of length five not is not spanned by faces. \end{defn} 

\begin{lem}\label{lem:noemptytrisqpent}A $\cat(0)$ simplicial $3$-complex contains no empty triangles, squares or pentagons. \end{lem}

\begin{proof} Let $\alpha = \{v_0, e_0, v_1, ..., v_{n} = v_0\}$ be an empty $n$-gon of minimal length,  $D$ be the minimal disk spanning $\alpha$. Since $\alpha$ is tight, $deg(v_i) >1$ for each $i$. If $deg(v_i) = 2$, then $v_i$ lies on a single triangle in $D$. Removing the triangle incident $v_i$ gives an $(n-1)$-gon $\alpha'$ spanned by $D' $ containing the same interior vertex as $D$, contradicting the choice of $\alpha$. Thus  $deg(v_i) \geq 3$ for all $i = 0, ..., n-1$, and $ n \geq \sum_{i=0}^{n-1} (4-deg(v_i)) \geq 6$. \end{proof} 

\begin{lem}\label{lem:fullnotempty} Full subcomplexes inherit the no empty triangles, squares and pentagons conditions.
\end{lem}

\begin{proof} Let $L$ be a full subcomplex of a simplicial complex $K$ with no empty triangles, squares or pentagons. The triangle condition follows immediately from the definition of full. Let $\alpha = [v_1,v_2,v_3,v_4,v_1]$ be a closed combinatorial path of length four in $L$. Then $\alpha$ is filled by a disk consisting of two triangular faces in $K$. Thus either $v_1$ and $v_3$, or $v_2$ and $v_4$ span an edge $e$ in $K$. The fullness condition implies $e$ is contained in $L$. Thus $\alpha$ must be spanned by two triangles sharing the edge $e$.  The same argument shows $L$ has no empty pentagons.
\end{proof}

\begin{thm} If $K$ is a $\cat(0)$ simplicial $3$-complex then $K$ is flag.
\end{thm}

\begin{proof}
Suppose that the vertices $v_1,v_2,...,v_n$ are pairwise connected by an edge in $K$. Lemma \ref{lem:noemptytrisqpent} implies that each triplet of vertices from the set $V = \{v_1,v_2,...v_n\}$ span a face. Let $W = \{w_1,w_2,w_3,w_4\}$ be a $4$-tuple of points of $V$. Then the full subcomplex of $K$ with vertex set $W$ contains the $2$-skeleton of a tetrahedron. As $\cat(0)$ spaces are contractible and contractions strictly reduce distance\cite{Bo92}, $W$ spans a tetrahedron of $K$. Let $U = \{u_1,u_2,u_3,u_4,u_5\} \subseteq V$ be distinct. Then the $1$-skeleton of the full subcomplex of $K$ with vertex set $U$ is a $K_5$ graph. Any $4$-tuple of points in $U$ span a tetrahedron. Let $e$ be an interior edge of the subcomplex. Then $e$ is contained in exactly three tetrahedra, and $lk_K(e)$ contains a loop $c$ of three edges corresponding the three tetrahedra. The length of each edge is the dihedral angle at $e$ in each tetrahedron, $\arccos(\frac{1}{3})$. Thus $\ell(c) = 3\arccos(\frac{1}{3}) < \frac{3\pi}{2}$, a contradiction to the $\cat(0)$ condition. This implies that there can be at most four distinct vertices pairwise connected by edges and they must span a simplex of $K$. Thus $K$ is flag. 
\end{proof}


\section{Links in Full Subcomplexes}\label{sec:links}


In this section is devoted to analyzing the structure of metric links in full subcomplexes.

\begin{prop}\label{prop:linksnotempty} Let $v$ be a vertex of a $\cat(0)$ simplicial $3$-complex $K$. Then $lk(v)$ has no empty triangles, squares or pentagons.
\end{prop}


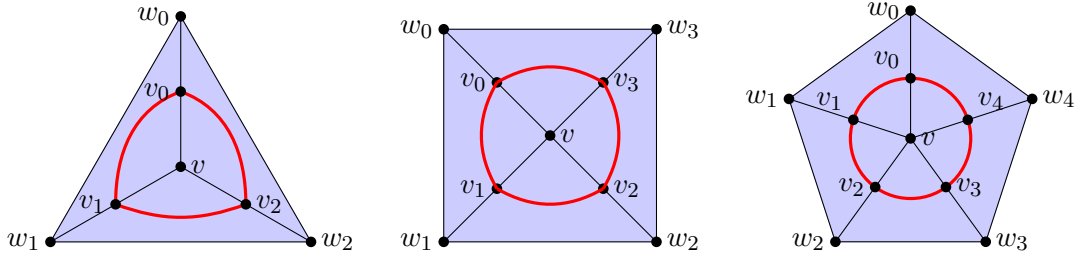
\begin{figure}

\begin{tabular}{ccc}

\begin{tikzpicture}

\coordinate (O) at (0:0);

\coordinate (A1) at (90: 2 cm);
\coordinate (A2) at (210: 2 cm);
\coordinate (A3) at (330: 2 cm);

\coordinate (B1) at (90: 1 cm);
\coordinate (B2) at (210: 1 cm);
\coordinate (B3) at (330: 1 cm);

\filldraw[blue!20] (A1)--(A2)--(A3)--cycle;

\draw (A1)--(A2)--(A3)--cycle;

\draw (O) -- (A1);
\draw (O)--(A2);
\draw (O)--(A3);

\draw[red, very thick] (B1) to [out =-160, in=90] (B2);
\draw[red, very thick] (B2) to [out= 340, in= 200] (B3);
\draw[red, very thick] (B3) to [out=90, in = -20] (B1);

\fill [black] (A1) circle (2pt)  node [anchor=east] {$w_0$};
\fill [black] (A2) circle (2pt) node [anchor=east] {$w_1$};
\fill [black] (A3) circle (2pt) node [anchor=west] {$w_2$};
\fill [black] (O) circle (2pt)  node [anchor=west] {$v$};

\fill [black] (B1) circle (2pt)  node [anchor=east] {$v_0$};
\fill [black] (B2) circle (2pt) node [anchor=east] {$v_1$};
\fill [black] (B3) circle (2pt) node [anchor=west] {$v_2$};

\end{tikzpicture} &

\begin{tikzpicture}

\coordinate (O) at (0,0);

\coordinate (A1) at (45*3: 2cm);
\coordinate (A2) at (45*5: 2cm);
\coordinate (A3) at (45*7: 2cm);
\coordinate (A4) at (45: 2cm);

\coordinate (B1) at (45*3: 1cm);
\coordinate (B2) at (45*5: 1cm);
\coordinate (B3) at (45*7: 1cm);
\coordinate (B4) at (45: 1cm);

\filldraw[blue!20] (A1)--(A2)--(A3)--(A4)--cycle;

\draw (A1)--(A2)--(A3)--(A4)--cycle; 

\draw (O) -- (A1);
\draw (O)--(A2);
\draw (O)--(A3);
\draw (O)--(A4);

\fill [black] (A1) circle (2pt)  node [anchor=east] {$w_0$};
\fill [black] (A2) circle (2pt) node [anchor=east] {$w_1$};
\fill [black] (A3) circle (2pt) node [anchor=west] {$w_2$};
\fill [black] (A4) circle (2pt) node [anchor=west] {$w_3$};
\fill [black] (O) circle (2pt)  node [anchor=west] {$v$};

\fill [black] (B1) circle (2pt)  node [anchor=east] {$v_0$};
\fill [black] (B2) circle (2pt) node [anchor=east] {$v_1$};
\fill [black] (B3) circle (2pt) node [anchor=west] {$v_2$};
\fill [black] (B4) circle (2pt) node [anchor=west] {$v_3$};

\draw[red,  very thick] (B1) to [out=240, in = 120] (B2);
\draw[red,  very thick] (B2) to [out = -30, in = 210] (B3);
\draw[red,  very thick] (B3) to [out = 60, in=300] (B4);
\draw[red,  very thick] (B4) to [out = 150, in = 30] (B1);

\end{tikzpicture} &

\begin{tikzpicture} 

\coordinate (A1) at (90: 1.7cm);
\coordinate (A2) at (90+72: 1.7cm);
\coordinate (A3) at (90+2*72: 1.7cm);
\coordinate (A4) at (90+3*72: 1.7cm);
\coordinate (A5) at (90+4*72: 1.7cm);

\coordinate (B1) at (90: 0.8cm);
\coordinate (B2) at (90+72: 0.8cm);
\coordinate (B3) at (90+2*72: 0.8cm);
\coordinate (B4) at (90+3*72: 0.8cm);
\coordinate (B5) at (90+4*72: 0.8cm);

\filldraw[blue!20] (A1)--(A2)--(A3)--(A4)--(A5)--cycle;

\draw (A1)--(A2)--(A3)--(A4)--(A5)--cycle; 

\draw (O) -- (A1);
\draw (O)--(A2);
\draw (O)--(A3);
\draw (O)--(A4);
\draw (O)--(A5);

\fill [black] (A1) circle (2pt)  node [anchor=east] {$w_0$};
\fill [black] (A2) circle (2pt) node [anchor=east] {$w_1$};
\fill [black] (A3) circle (2pt) node [anchor=east] {$w_2$};
\fill [black] (A4) circle (2pt) node [anchor=west] {$w_3$};
\fill [black] (A5) circle (2pt) node [anchor=west] {$w_4$};
\fill [black] (O) circle (2pt)  node [anchor=west] {$v$};

\draw[red,  very thick] (O) circle (0.8cm);

\fill [black] (B1) circle (2pt)  node [anchor=south east] {$v_0$};
\fill [black] (B2) circle (2pt) node [anchor= south east] {$v_1$};
\fill [black] (B3) circle (2pt) node [anchor=east] {$v_2$};
\fill [black] (B4) circle (2pt) node [anchor=west] {$v_3$};
\fill [black] (B5) circle (2pt) node [anchor=south west] {$v_4$};

\end{tikzpicture} \\

\end{tabular}

\caption{\label{fig:ngonlinks} Combinatorial triangles, squares and pentagons in metric links.}

\end{figure}


\begin{proof} Let $\alpha = [v_0,e_1,v_1,e_2,v_2,e_3,v_3=v_0] \subseteq lk(v)$ be a combinatorial triangle. Each vertex $v_i$ corresponds to an edge  in $K$ sharing  $v$ as a common vertex. Let $w_i$ be the vertex opposite $v$ on each edge.  The existence of $e_i $ in the link implies that $v$, $w_{i-1}$ and $w_{i}$ span a face of $K$, and $w_0$, $w_1$, $w_2$ lie on a combinatorial triangle in $K$ (see Figure \ref{fig:ngonlinks}). Thus the vertices $v$, $w_0$, $w_1$, and $w_2$ are pairwise connected by an edge and must span a simplex of $K$ by the flag condition. This implies $\{v_0,v_1,v_2\}$ spans a spherical triangle in $lk_K(v)$.

Next, let $\alpha$ be a combinatorial square in $lk(v)$ with vertices $v_0$, $v_1$, $v_2$, and $v_3$. As in the last case, the vertices $v_i$ correspond to edges in $K$ sharing $v$ as a vertex,  and the vertices $w_i$ opposing $v$  consecutively span edges. Thus $w_0$, $w_1$, $w_2$, and $w_3$ on lie on a combinatorial square in $K$. The complex $K$ has no empty squares, so either $\{w_0,w_2\}$ or $\{w_1,w_3\}$ spans an edge of $K$. Without loss of generality, assume $w_0$ and $w_2$ span an edge of $K$. Then $v$,$w_0$, $w_2$ are pairwise edge connected and span a face of $K$, implying $v_0$ and $v_2$ are connected by an edge in $lk(v)$. Thus $[v_0,v_1,v_2]$ and $[v_0,v_2,v_3]$ are combinatorial triangles in $lk(v)$, and $\{v_0,v_1,v_2\}$ and $\{v_0,v_2,v_3\}$ span spherical triangles.

It remains to show $lk(v)$ has no empty pentagons. Let $\alpha$ be a combinatorial pentagon in $lk(v)$ with vertices $v_0, v_1,..., v_4$. Again the vertices $v_i$ correspond to edges in $K$ sharing $v$,  and the vertices $w_0, w_1,...,w_4$ opposite $v$ form a combinatorial pentagon  in $K$. $K$ has no empty pentagons, so two pairs of non neighboring vertices on on the pentagon span edges. Without loss of generality suppose the sets $\{w_0,w_2\}$ and $\{w_0,w_3\}$ span edges. Then $\{v,w_0,w_2\}$ and $\{v,w_0,w_3\}$ span faces in $K$. Thus $\{v_0,v_2\}$ and $\{v_0,v_3\}$ span edges in $lk(v)$, implying  the sets $\{v_0,v_1,v_2\}$,$\{v_0,v_2,v_3\}$ and $\{v_0,v_3,v_4\}$ span faces in $lk(v)$. Thus $\alpha$ is nonempty.
\end{proof}

Combining Proposition \ref{prop:linksnotempty}, Lemma \ref{lem:linksfull} and  Lemma \ref{lem:fullnotempty} gives the following theorem.

\begin{thm}\label{prop:linksnotempty}Let $v$ be a vertex of a full subcomplex $L$ of a $\cat(0)$ simplicial complex. Then $lk_L(v)$ has no empty triangles, squares or pentagons.
\end{thm}


\section{Proof of Main Theorem}\label{sec:proof}


This section is devoted to the proof of the main theorem. By Gromov's link condition, nonpositive curvature is determined by the structure of the links of simplices of $L$. The link of a tetrahedron is empty and the link of a face is a discrete set of points. Neither contain short loops. Let $e$ be a edge of $L$. Then $lk_L(e)$ is a full subgraph of the metric graph $lk_K(e)$. Thus $lk_L(e)$ inherits the no short loop condition from $lk_K(e)$. Therefore for $3$-complexes, it suffices to analyze the structure of the links of vertices. This section begins with a brief discussion of curvature testing in metric complexes, and configurations. The details are found in \cite{ElMc04} and \cite{ElMc02}.

\begin{defn}[Configuration] A \emph{configuration} is a finite piecewise spherical complex $C$ which contains at least one locally geodesic loop $\rho$ of length less than $2\pi$.  A complex $K$ \emph{contains} a  configuration if there exist a cellular map $f: C \rightarrow lk_{K}(\sigma)$ for some simplex $\sigma \subseteq K$ such that $f(\rho)$ is also a locally geodesic loop of length less than $2\pi$ in $lk_{K}(\sigma)$. If $K$ does not contain a configuration $C$, then $K$ \emph{avoids} $C$. \end{defn}

In \cite{ElMc04}, Murray Elder and Jon McCammond prove that given a finite set of Euclidean polytopes, there exists a finite list of configurations that must be avoided in non positively curved complexes built out of the given polytopes. Precisely, Elder and McCammond show the following. 

\begin{thm}[Elder-McCammond] Let $K$ be a piecewise Euclidean complex with finitely many cell isometry types. Then there exists a finite list of configurations $\mathcal{C}$ such that $K$ is non positively curved if and if only $K$ avoids the configurations in $\mathcal{C}$. \end{thm}

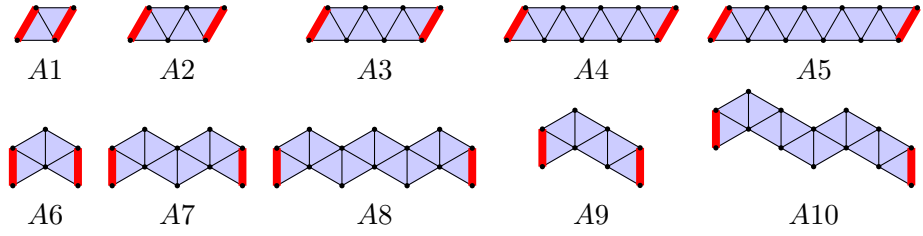
\begin{figure}[t]

\begin{center}

\begin{tabular}{ccccc}

\begin{tikzpicture}[scale=0.5]

\coordinate (A0) at (0,0);
\coordinate (A1) at (1,0);

\coordinate (B0) at (0.5, 0.5*{sqrt(3)});
\coordinate (B1) at (1.5, 0.5*{sqrt(3)});

\filldraw[blue!20] (A0)--(A1)--(B1)--(B0)--cycle;

\draw (A0)--(A1);
\draw (B0)--(B1);
\draw (B0)--(A1);

\draw[red,line width=3pt] (A0)--(B0);
\draw[red,line width=3pt] (A1)--(B1);

\fill [black] (A0) circle (2pt);
\fill [black] (A1) circle (2pt);
\fill [black] (B0) circle (2pt);
\fill [black] (B1) circle (2pt);

\end{tikzpicture}

&

\begin{tikzpicture}[scale=0.5]

\coordinate (A0) at (0,0);
\coordinate (A1) at (1,0);
\coordinate (A2) at (2,0);

\coordinate (B0) at (0.5, 0.5*{sqrt(3)});
\coordinate (B1) at (1.5, 0.5*{sqrt(3)});
\coordinate (B2) at (2.5, 0.5*{sqrt(3)});

\filldraw[blue!20] (A0)--(A2)--(B2)--(B0)--cycle;

\draw (A0)--(A2);
\draw (B0)--(B2);
\draw (B0)--(A1)--(B1)--(A2);

\draw[red,line width=3pt] (A0)--(B0);
\draw[red,line width=3pt] (A2)--(B2);

\fill [black] (A0) circle (2pt);
\fill [black] (A1) circle (2pt);
\fill [black] (A2) circle (2pt);
\fill [black] (B0) circle (2pt);
\fill [black] (B1) circle (2pt);
\fill [black] (B2) circle (2pt);

\end{tikzpicture}

&

\begin{tikzpicture}[scale=0.5]

\coordinate (A0) at (0,0);
\coordinate (A1) at (1,0);
\coordinate (A2) at (2,0);
\coordinate (A3) at (3,0);

\coordinate (B0) at (0.5, 0.5*{sqrt(3)});
\coordinate (B1) at (1.5, 0.5*{sqrt(3)});
\coordinate (B2) at (2.5, 0.5*{sqrt(3)});
\coordinate (B3) at (3.5, 0.5*{sqrt(3)});

\filldraw[blue!20] (A0)--(A3)--(B3)--(B0)--cycle;

\draw (A0)--(A3);
\draw (B0)--(B3);
\draw (B0)--(A1)--(B1)--(A2)--(B2)--(A3);

\draw[red,line width=3pt] (A0)--(B0);
\draw[red,line width=3pt] (A3)--(B3);

\fill [black] (A0) circle (2pt);
\fill [black] (A1) circle (2pt);
\fill [black] (A2) circle (2pt);
\fill [black] (A3) circle (2pt);

\fill [black] (B0) circle (2pt);
\fill [black] (B1) circle (2pt);
\fill [black] (B2) circle (2pt);
\fill [black] (B3) circle (2pt);

\end{tikzpicture}

&

\begin{tikzpicture}[scale=0.5]

\coordinate (A0) at (0,0);
\coordinate (A1) at (1,0);
\coordinate (A2) at (2,0);
\coordinate (A3) at (3,0);
\coordinate (A4) at (4,0);

\coordinate (B0) at (0.5, 0.5*{sqrt(3)});
\coordinate (B1) at (1.5, 0.5*{sqrt(3)});
\coordinate (B2) at (2.5, 0.5*{sqrt(3)});
\coordinate (B3) at (3.5, 0.5*{sqrt(3)});
\coordinate (B4) at (4.5, 0.5*{sqrt(3)});

\filldraw[blue!20] (A0)--(A4)--(B4)--(B0)--cycle;

\draw (A0)--(A4);
\draw (B0)--(B4);
\draw (B0)--(A1)--(B1)--(A2)--(B2)--(A3)--(B3)--(A4);

\draw[red,line width=3pt] (A0)--(B0);
\draw[red,line width=3pt] (A4)--(B4);

\fill [black] (A0) circle (2pt);
\fill [black] (A1) circle (2pt);
\fill [black] (A2) circle (2pt);
\fill [black] (A3) circle (2pt);
\fill [black] (A4) circle (2pt);

\fill [black] (B0) circle (2pt);
\fill [black] (B1) circle (2pt);
\fill [black] (B2) circle (2pt);
\fill [black] (B3) circle (2pt);
\fill [black] (B4) circle (2pt);

\end{tikzpicture}  &

\begin{tikzpicture}[scale=0.5]

\coordinate (A0) at (0,0);
\coordinate (A1) at (1,0);
\coordinate (A2) at (2,0);
\coordinate (A3) at (3,0);
\coordinate (A4) at (4,0);
\coordinate (A5) at (5,0);

\coordinate (B0) at (0.5, 0.5*{sqrt(3)});
\coordinate (B1) at (1.5, 0.5*{sqrt(3)});
\coordinate (B2) at (2.5, 0.5*{sqrt(3)});
\coordinate (B3) at (3.5, 0.5*{sqrt(3)});
\coordinate (B4) at (4.5, 0.5*{sqrt(3)});
\coordinate (B5) at (5.5, 0.5*{sqrt(3)});

\filldraw[blue!20] (A0)--(A5)--(B5)--(B0)--cycle;

\draw (A0)--(A5);
\draw (B0)--(B5);
\draw (B0)--(A1)--(B1)--(A2)--(B2)--(A3)--(B3)--(A4)--(B4)--(A5);

\draw[red,line width=3pt] (A0)--(B0);
\draw[red,line width=3pt] (A5)--(B5);

\fill [black] (A0) circle (2pt);
\fill [black] (A1) circle (2pt);
\fill [black] (A2) circle (2pt);
\fill [black] (A3) circle (2pt);
\fill [black] (A4) circle (2pt);
\fill [black] (A5) circle (2pt);

\fill [black] (B0) circle (2pt);
\fill [black] (B1) circle (2pt);
\fill [black] (B2) circle (2pt);
\fill [black] (B3) circle (2pt);
\fill [black] (B4) circle (2pt);
\fill [black] (B5) circle (2pt);

\end{tikzpicture} \\

$A1$ & $A2$ & $A3$ & $A4$ & $A5$\\

\begin{tikzpicture}[scale=0.5]

\coordinate (A0) at (0,0);
\coordinate (A2) at ({sqrt(3)},0);

\coordinate (B1) at (0.5*{sqrt(3)}, 0.5);

\coordinate (C0) at (0,1);
\coordinate (C2) at ({sqrt(3)},1);

\coordinate (D1) at (0.5*{sqrt(3)}, 1.5);

\filldraw[blue!20] (A0)--(B1)--(A2)--(C2)--(D1)--(C0)--cycle;

\draw (A0)--(B1)--(A2);
\draw (C0)--(D1)--(C2);
\draw (C0)--(B1)--(C2);
\draw (B1)--(D1);

\draw[red,line width=3pt] (A0)--(C0);
\draw[red,line width=3pt] (A2)--(C2);

\fill [black] (A0) circle (2pt);
\fill [black] (A2) circle (2pt);

\fill [black] (C0) circle (2pt);
\fill [black] (C2) circle (2pt);

\fill [black] (B1) circle (2pt);

\fill [black] (D1) circle (2pt);

\end{tikzpicture}

&

\begin{tikzpicture}[scale=0.5]

\coordinate (A0) at (0,0);
\coordinate (A2) at ({sqrt(3)},0);
\coordinate (A4) at (2*{sqrt(3)},0);

\coordinate (B1) at (0.5*{sqrt(3)}, 0.5);
\coordinate (B3) at (1.5*{sqrt(3)}, 0.5);

\coordinate (C0) at (0,1);
\coordinate (C2) at ({sqrt(3)},1);
\coordinate (C4) at (2*{sqrt(3)},1);

\coordinate (D1) at (0.5*{sqrt(3)}, 1.5);
\coordinate (D3) at (1.5*{sqrt(3)}, 1.5);

\filldraw[blue!20] (A0)--(B1)--(A2)--(B3)--(A4)--(C4)--(D3)--(C2)--(D1)--(C0)--cycle;

\draw (A0)--(B1)--(A2)--(B3)--(A4);
\draw (C0)--(D1)--(C2)--(D3)--(C4);
\draw (C0)--(B1)--(C2)--(B3)--(C4);
\draw (B1)--(D1);
\draw (A2)--(C2);
\draw (B3)--(D3);

\draw[red,line width=3pt] (A0)--(C0);
\draw[red,line width=3pt] (A4)--(C4);

\fill [black] (A0) circle (2pt);
\fill [black] (A2) circle (2pt);
\fill [black] (A4) circle (2pt);

\fill [black] (C0) circle (2pt);
\fill [black] (C2) circle (2pt);
\fill [black] (C4) circle (2pt);

\fill [black] (B1) circle (2pt);
\fill [black] (B3) circle (2pt);

\fill [black] (D1) circle (2pt);
\fill [black] (D3) circle (2pt);

\end{tikzpicture}

&

\begin{tikzpicture}[scale=0.5]

\coordinate (A0) at (0,0);
\coordinate (A2) at ({sqrt(3)},0);
\coordinate (A4) at (2*{sqrt(3)},0);
\coordinate (A6) at (3*{sqrt(3)},0);

\coordinate (B1) at (0.5*{sqrt(3)}, 0.5);
\coordinate (B3) at (1.5*{sqrt(3)}, 0.5);
\coordinate (B5) at (2.5*{sqrt(3)}, 0.5);

\coordinate (C0) at (0,1);
\coordinate (C2) at ({sqrt(3)},1);
\coordinate (C4) at (2*{sqrt(3)},1);
\coordinate (C6) at (3*{sqrt(3)},1);

\coordinate (D1) at (0.5*{sqrt(3)}, 1.5);
\coordinate (D3) at (1.5*{sqrt(3)}, 1.5);
\coordinate (D5) at (2.5*{sqrt(3)}, 1.5);

\filldraw[blue!20] (A0)--(B1)--(A2)--(B3)--(A4)--(B5)--(A6)--(C6)--(D5)--(C4)--(D3)--(C2)--(D1)--(C0)--cycle;

\draw (A0)--(B1)--(A2)--(B3)--(A4)--(B5)--(A6);
\draw (C0)--(D1)--(C2)--(D3)--(C4)--(D5)--(C6);
\draw (C0)--(B1)--(C2)--(B3)--(C4)--(B5)--(C6);
\draw (B1)--(D1);
\draw (A2)--(C2);
\draw (B3)--(D3);
\draw (A4)--(C4);
\draw (B5)--(D5);

\draw[red,line width=3pt] (A0)--(C0);
\draw[red,line width=3pt] (A6)--(C6);

\fill [black] (A0) circle (2pt);
\fill [black] (A2) circle (2pt);
\fill [black] (A4) circle (2pt);
\fill [black] (A6) circle (2pt);

\fill [black] (C0) circle (2pt);
\fill [black] (C2) circle (2pt);
\fill [black] (C4) circle (2pt);
\fill [black] (C6) circle (2pt);

\fill [black] (B1) circle (2pt);
\fill [black] (B3) circle (2pt);
\fill [black] (B5) circle (2pt);

\fill [black] (D1) circle (2pt);
\fill [black] (D3) circle (2pt);
\fill [black] (D5) circle (2pt);

\end{tikzpicture} &

\begin{tikzpicture}[scale=0.5]

\coordinate (A0) at (0,0);
\coordinate (A2) at ({sqrt(3)},0);

\coordinate (B1) at (0.5*{sqrt(3)}, 0.5);
\coordinate (B3) at (1.5*{sqrt(3)}, 0.5);

\coordinate (C0) at (0,1);
\coordinate (C2) at ({sqrt(3)},1);

\coordinate (D1) at (0.5*{sqrt(3)}, 1.5);

\coordinate (L3) at (1.5*{sqrt(3)}, -0.5);

\filldraw[blue!20] (A0)--(B1)--(A2)--(L3)--(B3)--(C2)--(D1)--(C0)--cycle;

\draw (A0)--(B1)--(A2)--(L3);
\draw (C0)--(D1)--(C2)--(B3);
\draw (C0)--(B1)--(C2);
\draw (A2)--(B3);
\draw (B1)--(D1);
\draw (A2)--(C2);

\draw[red,line width=3pt] (A0)--(C0);
\draw[red,line width=3pt] (L3)--(B3);

\fill [black] (A0) circle (2pt);
\fill [black] (A2) circle (2pt);

\fill [black] (C0) circle (2pt);
\fill [black] (C2) circle (2pt);

\fill [black] (B1) circle (2pt);
\fill [black] (B3) circle (2pt);

\fill [black] (D1) circle (2pt);

\fill [black] (L3) circle (2pt);

\end{tikzpicture}

&

\begin{tikzpicture}[scale=0.5]

\coordinate (A0) at (0,0);
\coordinate (A2) at ({sqrt(3)},0);
\coordinate (A4) at (2*{sqrt(3)},0);
\coordinate (A6) at (3*{sqrt(3)},0);

\coordinate (B1) at (0.5*{sqrt(3)}, 0.5);
\coordinate (B3) at (1.5*{sqrt(3)}, 0.5);
\coordinate (B5) at (2.5*{sqrt(3)}, 0.5);

\coordinate (C0) at (0,1);
\coordinate (C2) at ({sqrt(3)},1);
\coordinate (C4) at (2*{sqrt(3)},1);

\coordinate (D1) at (0.5*{sqrt(3)}, 1.5);

\coordinate (L3) at (1.5*{sqrt(3)}, -0.5);
\coordinate (L5) at (2.5*{sqrt(3)}, -0.5);

\coordinate (M6) at (3*{sqrt(3)},-1);

\filldraw[blue!20] (A0)--(B1)--(A2)--(L3)--(A4)--(L5)--(M6)--(A6)--(B5)--(C4)--(B3)--(C2)--(D1)--(C0)--cycle;

\draw (A0)--(B1)--(A2)--(L3)--(A4)--(L5)--(M6);
\draw (C0)--(D1)--(C2)--(B3)--(C4)--(B5)--(A6);
\draw (C0)--(B1)--(C2);
\draw (A2)--(B3)--(A4)--(B5);
\draw (L5)--(A6);
\draw (B1)--(D1);
\draw (A2)--(C2);
\draw (L3)--(B3);
\draw (A4)--(C4);
\draw (L5)--(B5);

\draw[red,line width=3pt] (A0)--(C0);
\draw[red,line width=3pt] (M6)--(A6);

\fill [black] (A0) circle (2pt);
\fill [black] (A2) circle (2pt);
\fill [black] (A4) circle (2pt);
\fill [black] (A6) circle (2pt);

\fill [black] (C0) circle (2pt);
\fill [black] (C2) circle (2pt);
\fill [black] (C4) circle (2pt);

\fill [black] (B1) circle (2pt);
\fill [black] (B3) circle (2pt);
\fill [black] (B5) circle (2pt);

\fill [black] (D1) circle (2pt);

\fill [black] (L3) circle (2pt);
\fill [black] (L5) circle (2pt);

\fill [black] (M6) circle (2pt);

\end{tikzpicture}  \\

 $A6$ & $A7$ & $A8$ & $A9$ & $A10$ \\
 
 \end{tabular}

\caption{\label{fig:annular} The annular galleries described in \cite{ElMc02}. Annular galleries are built by identifying the thick, red outer edges.} 

\end{center}

\end{figure}

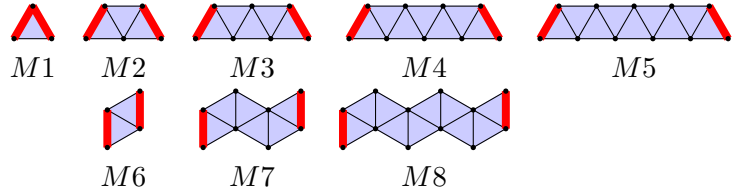
\begin{figure}[t]

\begin{center}

\begin{tabular}{ccccc}

\begin{tikzpicture}[scale=0.5]

\coordinate (A0) at (0,0);
\coordinate (A1) at (1,0);

\coordinate (B0) at (0.5, 0.5*{sqrt(3)});

\filldraw[blue!20] (A0)--(A1)--(B0)--cycle;

\draw (A0)--(A1);

\draw [red, line width=3pt] (A0)--(B0);
\draw [red, line width=3pt] (A1)--(B0);

\fill [black] (A0) circle (2pt);
\fill [black] (A1) circle (2pt);

\fill [black] (B0) circle (2pt);

\end{tikzpicture} &
 
 \begin{tikzpicture}[scale=0.5]

\coordinate (A0) at (0,0);
\coordinate (A1) at (1,0);
\coordinate (A2) at (2,0);

\coordinate (B0) at (0.5, 0.5*{sqrt(3)});
\coordinate (B1) at (1.5, 0.5*{sqrt(3)});

\filldraw[blue!20] (A0)--(A1)--(A2)--(B1)--(B0)--cycle;

\draw (A0)--(A2);
\draw (B0)--(B1);
\draw (B0)--(A1)--(B1);

\draw [red, line width=3pt] (A0)--(B0);
\draw [red, line width=3pt] (A2)--(B1);

\fill [black] (A0) circle (2pt);
\fill [black] (A1) circle (2pt);
\fill [black] (A2) circle (2pt);

\fill [black] (B0) circle (2pt);
\fill [black] (B1) circle (2pt);

\end{tikzpicture} &

\begin{tikzpicture}[scale=0.5]

\coordinate (A0) at (0,0);
\coordinate (A1) at (1,0);
\coordinate (A2) at (2,0);
\coordinate (A3) at (3,0);

\coordinate (B0) at (0.5, 0.5*{sqrt(3)});
\coordinate (B1) at (1.5, 0.5*{sqrt(3)});
\coordinate (B2) at (2.5, 0.5*{sqrt(3)});

\filldraw[blue!20] (A0)--(A1)--(A2)--(A3)--(B2)--(B1)--(B0)--cycle;

\draw (A0)--(A3);
\draw (B0)--(B2);
\draw (B0)--(A1)--(B1)--(A2)--(B2);

\draw [red, line width=3pt] (A0)--(B0);
\draw [red, line width=3pt] (A3)--(B2);

\fill [black] (A0) circle (2pt);
\fill [black] (A1) circle (2pt);
\fill [black] (A2) circle (2pt);
\fill [black] (A3) circle (2pt);

\fill [black] (B0) circle (2pt);
\fill [black] (B1) circle (2pt);
\fill [black] (B2) circle (2pt);

\end{tikzpicture} &

\begin{tikzpicture}[scale=0.5]

\coordinate (A0) at (0,0);
\coordinate (A1) at (1,0);
\coordinate (A2) at (2,0);
\coordinate (A3) at (3,0);
\coordinate (A4) at (4,0);

\coordinate (B0) at (0.5, 0.5*{sqrt(3)});
\coordinate (B1) at (1.5, 0.5*{sqrt(3)});
\coordinate (B2) at (2.5, 0.5*{sqrt(3)});
\coordinate (B3) at (3.5, 0.5*{sqrt(3)});

\filldraw[blue!20] (A0)--(A1)--(A2)--(A3)--(A4)--(B3)--(B2)--(B1)--(B0)--cycle;

\draw (A0)--(A4);
\draw (B0)--(B3);
\draw (B0)--(A1)--(B1)--(A2)--(B2)--(A3)--(B3);

\draw [red, line width=3pt] (A0)--(B0);
\draw [red, line width=3pt] (A4)--(B3);

\fill [black] (A0) circle (2pt);
\fill [black] (A1) circle (2pt);
\fill [black] (A2) circle (2pt);
\fill [black] (A3) circle (2pt);
\fill [black] (A4) circle (2pt);

\fill [black] (B0) circle (2pt);
\fill [black] (B1) circle (2pt);
\fill [black] (B2) circle (2pt);
\fill [black] (B3) circle (2pt);

\end{tikzpicture}  &

\begin{tikzpicture}[scale=0.5]

\coordinate (A0) at (0,0);
\coordinate (A1) at (1,0);
\coordinate (A2) at (2,0);
\coordinate (A3) at (3,0);
\coordinate (A4) at (4,0);
\coordinate (A5) at (5,0);

\coordinate (B0) at (0.5, 0.5*{sqrt(3)});
\coordinate (B1) at (1.5, 0.5*{sqrt(3)});
\coordinate (B2) at (2.5, 0.5*{sqrt(3)});
\coordinate (B3) at (3.5, 0.5*{sqrt(3)});
\coordinate (B4) at (4.5, 0.5*{sqrt(3)});

\filldraw[blue!20] (A0)--(A1)--(A2)--(A3)--(A4)--(A5)--(B4)--(B3)--(B2)--(B1)--(B0)--cycle;

\draw (A0)--(A5);
\draw (B0)--(B4);
\draw (B0)--(A1)--(B1)--(A2)--(B2)--(A3)--(B3)--(A4)--(B4);

\draw [red, line width=3pt] (A0)--(B0);
\draw [red, line width=3pt] (A5)--(B4);

\fill [black] (A0) circle (2pt);
\fill [black] (A1) circle (2pt);
\fill [black] (A2) circle (2pt);
\fill [black] (A3) circle (2pt);
\fill [black] (A4) circle (2pt);
\fill [black] (A5) circle (2pt);

\fill [black] (B0) circle (2pt);
\fill [black] (B1) circle (2pt);
\fill [black] (B2) circle (2pt);
\fill [black] (B3) circle (2pt);
\fill [black] (B4) circle (2pt);

\end{tikzpicture}  \\

$M1$ & $M2$ & $M3$ & $M4$ & $M5$ \\

&

\begin{tikzpicture}[scale=0.5]

\coordinate (A0) at (0,0);

\coordinate (B1) at (0.5*{sqrt(3)}, 0.5);

\coordinate (C0) at (0,1);

\coordinate (D1) at (0.5*{sqrt(3)}, 1.5);

\filldraw[blue!20] (A0)--(B1)--(D1)--(C0)--cycle;

\draw (A0)--(B1);
\draw (C0)--(D1);
\draw (C0)--(B1);
\draw (B1)--(D1);

\draw[red,line width=3pt] (A0)--(C0);
\draw[red,line width=3pt] (B1)--(D1);

\fill [black] (A0) circle (2pt);

\fill [black] (C0) circle (2pt);

\fill [black] (B1) circle (2pt);

\fill [black] (D1) circle (2pt);

\end{tikzpicture} &

\begin{tikzpicture}[scale=0.5]

\coordinate (A0) at (0,0);
\coordinate (A2) at ({sqrt(3)},0);

\coordinate (B1) at (0.5*{sqrt(3)}, 0.5);
\coordinate (B3) at (1.5*{sqrt(3)}, 0.5);

\coordinate (C0) at (0,1);
\coordinate (C2) at ({sqrt(3)},1);

\coordinate (D1) at (0.5*{sqrt(3)}, 1.5);
\coordinate (D3) at (1.5*{sqrt(3)}, 1.5);

\filldraw[blue!20] (A0)--(B1)--(A2)--(B3)--(D3)--(C2)--(D1)--(C0)--cycle;

\draw (A0)--(B1)--(A2)--(B3);
\draw (C0)--(D1)--(C2)--(D3);
\draw (C0)--(B1)--(C2)--(B3);
\draw (B1)--(D1);
\draw (A2)--(C2);

\draw[red,line width=3pt] (A0)--(C0);
\draw[red,line width=3pt] (B3)--(D3);

\fill [black] (A0) circle (2pt);
\fill [black] (A2) circle (2pt);

\fill [black] (C0) circle (2pt);
\fill [black] (C2) circle (2pt);

\fill [black] (B1) circle (2pt);
\fill [black] (B3) circle (2pt);

\fill [black] (D1) circle (2pt);
\fill [black] (D3) circle (2pt);

\end{tikzpicture} &

\begin{tikzpicture}[scale=0.5]

\coordinate (A0) at (0,0);
\coordinate (A2) at ({sqrt(3)},0);
\coordinate (A4) at (2*{sqrt(3)},0);

\coordinate (B1) at (0.5*{sqrt(3)}, 0.5);
\coordinate (B3) at (1.5*{sqrt(3)}, 0.5);
\coordinate (B5) at (2.5*{sqrt(3)}, 0.5);

\coordinate (C0) at (0,1);
\coordinate (C2) at ({sqrt(3)},1);
\coordinate (C4) at (2*{sqrt(3)},1);

\coordinate (D1) at (0.5*{sqrt(3)}, 1.5);
\coordinate (D3) at (1.5*{sqrt(3)}, 1.5);
\coordinate (D5) at (2.5*{sqrt(3)}, 1.5);

\filldraw[blue!20] (A0)--(B1)--(A2)--(B3)--(A4)--(B5)--(D5)--(C4)--(D3)--(C2)--(D1)--(C0)--cycle;

\draw (A0)--(B1)--(A2)--(B3)--(A4)--(B5);
\draw (C0)--(D1)--(C2)--(D3)--(C4)--(D5);
\draw (C0)--(B1)--(C2)--(B3)--(C4)--(B5);
\draw (B1)--(D1);
\draw (A2)--(C2);
\draw (B3)--(D3);
\draw (A4)--(C4);

\draw[red,line width=3pt] (A0)--(C0);
\draw[red,line width=3pt] (B5)--(D5);

\fill [black] (A0) circle (2pt);
\fill [black] (A2) circle (2pt);
\fill [black] (A4) circle (2pt);

\fill [black] (C0) circle (2pt);
\fill [black] (C2) circle (2pt);
\fill [black] (C4) circle (2pt);

\fill [black] (B1) circle (2pt);
\fill [black] (B3) circle (2pt);
\fill [black] (B5) circle (2pt);

\fill [black] (D1) circle (2pt);
\fill [black] (D3) circle (2pt);
\fill [black] (D5) circle (2pt);

\end{tikzpicture} \\

 & $ M6$ & $M7$ & $M8$ &  \\

\end{tabular}

\end{center}

\caption{\label{fig:mobius} The m\"{o}bius galleries described by Elder and McCammond. To build m\"{o}bius galleries outer edges are again identified but with a twist.}

\end{figure}

\begin{figure}[t]

\begin{center}

\begin{tabular}{ccccc}

\begin{tikzpicture}[scale=0.5]

\draw (0,0)--(1,0);

\fill [black] (0,0) circle (2pt);
\fill [black] (1,0) circle (2pt);

\end{tikzpicture} &

\begin{tikzpicture}[scale=0.5]

\draw (0,0)--(2,0);

\fill [black] (0,0) circle (2pt);
\fill [black] (1,0) circle (2pt);
\fill [black] (2,0) circle (2pt);

\end{tikzpicture} & 

\begin{tikzpicture}[scale=0.5]

\draw (0,0)--(3,0);

\fill [black] (0,0) circle (2pt);
\fill [black] (1,0) circle (2pt);
\fill [black] (2,0) circle (2pt);
\fill [black] (3,0) circle (2pt);

\end{tikzpicture} &

\begin{tikzpicture}[scale=0.5]

\draw (0,0)--(4,0);

\fill [black] (0,0) circle (2pt);
\fill [black] (1,0) circle (2pt);
\fill [black] (2,0) circle (2pt);
\fill [black] (3,0) circle (2pt);
\fill [black] (4,0) circle (2pt);

\end{tikzpicture} &

\begin{tikzpicture}[scale=0.5]

\draw (0,0)--(5,0);

\fill [black] (0,0) circle (2pt);
\fill [black] (1,0) circle (2pt);
\fill [black] (2,0) circle (2pt);
\fill [black] (3,0) circle (2pt);
\fill [black] (4,0) circle (2pt);
\fill [black] (5,0) circle (2pt);

\end{tikzpicture} \\

$N1$ & $N2$ & $N3$ & $N4$ & $N5$ \\

\begin{tikzpicture}[scale=0.5]

\coordinate (A0) at (0,0);
\coordinate (A2) at ({sqrt(3)},0);

\coordinate (B1) at (0.5*{sqrt(3)}, 0.5);

\coordinate (L1) at (0.5*{sqrt(3)}, -0.5);

\filldraw[blue!20] (A0)--(L1)--(A2)--(B1)--cycle;

\draw (A0)--(L1)--(A2)--(B1)--cycle;

\draw (L1)--(B1);

\fill [red] (A0) circle (3pt);
\fill [red] (A2) circle (3pt);

\fill [black] (B1) circle (2pt);
\fill [black] (L1) circle (2pt);

\end{tikzpicture} &

\begin{tikzpicture}[scale=0.5]

\coordinate (A0) at (0,0);
\coordinate (A2) at ({sqrt(3)},0);
\coordinate (A3) at ({sqrt(3)}+1,0);

\coordinate (B1) at (0.5*{sqrt(3)}, 0.5);

\coordinate (L1) at (0.5*{sqrt(3)}, -0.5);

\filldraw[blue!20] (A0)--(L1)--(A2)--(B1)--cycle;

\draw (A0)--(L1)--(A2)--(B1)--cycle;
\draw (A2)--(A3);

\draw (L1)--(B1);

\fill [red] (A0) circle (3pt);
\fill [red] (A3) circle (3pt);

\fill [black] (A2) circle (2pt);

\fill [black] (B1) circle (2pt);
\fill [black] (L1) circle (2pt);

\end{tikzpicture} &

\begin{tikzpicture}[scale=0.5]

\coordinate (A0) at (0,0);
\coordinate (A2) at ({sqrt(3)},0);
\coordinate (A3) at ({sqrt(3)}+1,0);
\coordinate (A4) at ({sqrt(3)}+2,0);

\coordinate (B1) at (0.5*{sqrt(3)}, 0.5);

\coordinate (L1) at (0.5*{sqrt(3)}, -0.5);

\filldraw[blue!20] (A0)--(L1)--(A2)--(B1)--cycle;

\draw (A0)--(L1)--(A2)--(B1)--cycle;
\draw (A2)--(A4);

\draw (L1)--(B1);

\fill [red] (A0) circle (3pt);
\fill [red] (A4) circle (3pt);

\fill [black] (A2) circle (2pt);
\fill [black] (A3) circle (2pt);

\fill [black] (B1) circle (2pt);
\fill [black] (L1) circle (2pt);

\end{tikzpicture} &

\begin{tikzpicture}[scale=0.5]

\coordinate (A0) at (0,0);
\coordinate (A2) at ({sqrt(3)},0);
\coordinate (A3) at ({sqrt(3)}+1,0);
\coordinate (A4) at ({sqrt(3)}+2,0);
\coordinate (A5) at ({sqrt(3)}+3,0);

\coordinate (B1) at (0.5*{sqrt(3)}, 0.5);

\coordinate (L1) at (0.5*{sqrt(3)}, -0.5);

\filldraw[blue!20] (A0)--(L1)--(A2)--(B1)--cycle;

\draw (A0)--(L1)--(A2)--(B1)--cycle;
\draw (A2)--(A5);

\draw (L1)--(B1);

\fill [red] (A0) circle (3pt);
\fill [red] (A5) circle (3pt);

\fill [black] (A2) circle (2pt);
\fill [black] (A3) circle (2pt);
\fill [black] (A4) circle (2pt);

\fill [black] (B1) circle (2pt);
\fill [black] (L1) circle (2pt);

\end{tikzpicture} &

\begin{tikzpicture}[scale=0.5]

\coordinate (A0) at (0,0);
\coordinate (A2) at ({sqrt(3)},0);
\coordinate (A3) at ({sqrt(3)}+1,0);
\coordinate (A4) at ({sqrt(3)}+2,0);
\coordinate (A5) at ({sqrt(3)}+3,0);
\coordinate (A6) at ({sqrt(3)}+4,0);

\coordinate (B1) at (0.5*{sqrt(3)}, 0.5);

\coordinate (L1) at (0.5*{sqrt(3)}, -0.5);

\filldraw[blue!20] (A0)--(L1)--(A2)--(B1)--cycle;

\draw (A0)--(L1)--(A2)--(B1)--cycle;
\draw (A2)--(A6);

\draw (L1)--(B1);

\fill [red] (A0) circle (3pt);
\fill [red] (A6) circle (3pt);

\fill [black] (A2) circle (2pt);
\fill [black] (A3) circle (2pt);
\fill [black] (A4) circle (2pt);
\fill [black] (A5) circle (2pt);

\fill [black] (B1) circle (2pt);
\fill [black] (L1) circle (2pt);

\end{tikzpicture} \\

$N6$ & $N7$ & $N8$ & $N9$ & $N10$\\

\begin{tikzpicture}[scale=0.5]

\coordinate (A0) at (0,0);
\coordinate (A2) at ({sqrt(3)},0);
\coordinate (A4) at (2*{sqrt(3)},0);

\coordinate (B1) at (0.5*{sqrt(3)}, 0.5);
\coordinate (B3) at (1.5*{sqrt(3)}, 0.5);

\coordinate (L1) at (0.5*{sqrt(3)}, -0.5);
\coordinate (L3) at (1.5*{sqrt(3)}, -0.5);

\filldraw[blue!20] (A0)--(L1)--(A2)--(L3)--(A4)--(B3)--(A2)--(B1)--cycle;

\draw (A0)--(L1)--(A2)--(L3)--(A4)--(B3)--(A2)--(B1)--cycle;

\draw (L1)--(B1);
\draw (L3)--(B3);

\fill [red] (A0) circle (3pt);
\fill [red] (A4) circle (3pt);

\fill [black] (A2) circle (2pt);

\fill [black] (B1) circle (2pt);
\fill [black] (B3) circle (2pt);

\fill [black] (L1) circle (2pt);
\fill [black] (L3) circle (2pt);

\end{tikzpicture} &

\begin{tikzpicture}[scale=0.5]

\coordinate (A0) at (0,0);
\coordinate (A2) at ({sqrt(3)},0);
\coordinate (A4) at (2*{sqrt(3)},0);
\coordinate (A5) at (2*{sqrt(3)}+1,0);

\coordinate (B1) at (0.5*{sqrt(3)}, 0.5);
\coordinate (B3) at (1.5*{sqrt(3)}, 0.5);

\coordinate (L1) at (0.5*{sqrt(3)}, -0.5);
\coordinate (L3) at (1.5*{sqrt(3)}, -0.5);

\filldraw[blue!20] (A0)--(L1)--(A2)--(L3)--(A4)--(B3)--(A2)--(B1)--cycle;

\draw (A0)--(L1)--(A2)--(L3)--(A4)--(B3)--(A2)--(B1)--cycle;

\draw (L1)--(B1);
\draw (L3)--(B3);

\draw (A4)--(A5);

\fill [red] (A0) circle (3pt);
\fill [red] (A5) circle (3pt);

\fill [black] (A2) circle (2pt);
\fill [black] (A4) circle (2pt);

\fill [black] (B1) circle (2pt);
\fill [black] (B3) circle (2pt);

\fill [black] (L1) circle (2pt);
\fill [black] (L3) circle (2pt);

\end{tikzpicture} &

\begin{tikzpicture}[scale=0.5]

\coordinate (A0) at (0,0);
\coordinate (A2) at ({sqrt(3)},0);
\coordinate (A4) at (2*{sqrt(3)},0);
\coordinate (A5) at (2*{sqrt(3)}+1,0);
\coordinate (A6) at (2*{sqrt(3)}+2,0);

\coordinate (B1) at (0.5*{sqrt(3)}, 0.5);
\coordinate (B3) at (1.5*{sqrt(3)}, 0.5);

\coordinate (L1) at (0.5*{sqrt(3)}, -0.5);
\coordinate (L3) at (1.5*{sqrt(3)}, -0.5);

\filldraw[blue!20] (A0)--(L1)--(A2)--(L3)--(A4)--(B3)--(A2)--(B1)--cycle;

\draw (A0)--(L1)--(A2)--(L3)--(A4)--(B3)--(A2)--(B1)--cycle;

\draw (L1)--(B1);
\draw (L3)--(B3);

\draw (A4)--(A6);

\fill [red] (A0) circle (3pt);
\fill [red] (A6) circle (3pt);

\fill [black] (A2) circle (2pt);
\fill [black] (A4) circle (2pt);
\fill [black] (A5) circle (2pt);

\fill [black] (B1) circle (2pt);
\fill [black] (B3) circle (2pt);

\fill [black] (L1) circle (2pt);
\fill [black] (L3) circle (2pt);

\end{tikzpicture} &

\begin{tikzpicture}[scale=0.5]

\coordinate (A0) at (0,0);
\coordinate (A2) at ({sqrt(3)},0);
\coordinate (A3) at ({sqrt(3)}+1,0);
\coordinate (A5) at (2*{sqrt(3)}+1,0);
\coordinate (A6) at (2*{sqrt(3)}+2,0);

\coordinate (B1) at (0.5*{sqrt(3)}, 0.5);
\coordinate (B4) at (1.5*{sqrt(3)}+1, 0.5);

\coordinate (L1) at (0.5*{sqrt(3)}, -0.5);
\coordinate (L4) at (1.5*{sqrt(3)}+1, -0.5);

\filldraw[blue!20] (A0)--(L1)--(A2)--(B1)--cycle;
\filldraw[blue!20] (A3)--(L4)--(A5)--(B4)--cycle;

\draw (A0)--(L1)--(A2)--(B1)--cycle;
\draw (A3)--(L4)--(A5)--(B4)--cycle;

\draw (L1)--(B1);
\draw (L4)--(B4);

\draw (A2)--(A3);
\draw (A5)--(A6);

\fill [red] (A0) circle (3pt);
\fill [red] (A6) circle (3pt);

\fill [black] (A2) circle (2pt);
\fill [black] (A3) circle (2pt);
\fill [black] (A5) circle (2pt);

\fill [black] (B1) circle (2pt);
\fill [black] (B4) circle (2pt);

\fill [black] (L1) circle (2pt);
\fill [black] (L4) circle (2pt);

\end{tikzpicture} &

\begin{tikzpicture}[scale=0.5]

\coordinate (A0) at (0,0);
\coordinate (A2) at ({sqrt(3)},0);
\coordinate (A4) at (2*{sqrt(3)},0);
\coordinate (A6) at (3*{sqrt(3)},0);

\coordinate (B1) at (0.5*{sqrt(3)}, 0.5);
\coordinate (B3) at (1.5*{sqrt(3)}, 0.5);
\coordinate (B5) at (2.5*{sqrt(3)}, 0.5);

\coordinate (L1) at (0.5*{sqrt(3)}, -0.5);
\coordinate (L3) at (1.5*{sqrt(3)}, -0.5);
\coordinate (L5) at (2.5*{sqrt(3)}, -0.5);

\filldraw[blue!20] (A0)--(L1)--(A2)--(L3)--(A4)--(L5)--(A6)--(B5)--(A4)--(B3)--(A2)--(B1)--cycle;

\draw (A0)--(L1)--(A2)--(L3)--(A4)--(L5)--(A6)--(B5)--(A4)--(B3)--(A2)--(B1)--cycle;

\draw (L1)--(B1);
\draw (L3)--(B3);
\draw (L5)--(B5);

\fill [red] (A0) circle (3pt);
\fill [red] (A6) circle (3pt);

\fill [black] (A2) circle (2pt);
\fill [black] (A4) circle (2pt);

\fill [black] (B1) circle (2pt);
\fill [black] (B3) circle (2pt);
\fill [black] (B5) circle (2pt);

\fill [black] (L1) circle (2pt);
\fill [black] (L3) circle (2pt);
\fill [black] (L5) circle (2pt);

\end{tikzpicture} \\
$N11$ & $N12$ & $N13$ & $N14$ & $N15$\\

\begin{tikzpicture}[scale=0.5]

\coordinate (A0) at (0,0);
\coordinate (A2) at ({sqrt(3)},0);

\coordinate (B1) at (0.5*{sqrt(3)}, 0.5);
\coordinate (B3) at (1.5*{sqrt(3)}, 0.5);

\coordinate (C2) at ({sqrt(3)},1);

\coordinate (L1) at (0.5*{sqrt(3)}, -0.5);

\filldraw[blue!20] (A0)--(L1)--(B3)--(C2)--cycle;

\draw (A0)--(L1)--(B3)--(C2)--cycle;
\draw (L1)--(B1)--(A2)--(C2);

\fill [red] (A0) circle (3pt);
\fill [red] (B3) circle (3pt);

\fill [black] (A2) circle (2pt);

\fill [black] (B1) circle (2pt);

\fill [black] (C2) circle (2pt);

\fill [black] (L1) circle (2pt);

\end{tikzpicture} &

\begin{tikzpicture}[scale=0.5]

\coordinate (A0) at (0,0);
\coordinate (A2) at ({sqrt(3)},0);

\coordinate (B1) at (0.5*{sqrt(3)}, 0.5);
\coordinate (B3) at (1.5*{sqrt(3)}, 0.5);
\coordinate (B4) at (1.5*{sqrt(3)}+1, 0.5);

\coordinate (C2) at ({sqrt(3)},1);

\coordinate (L1) at (0.5*{sqrt(3)}, -0.5);

\filldraw[blue!20] (A0)--(L1)--(B3)--(C2)--cycle;

\draw (A0)--(L1)--(B3)--(C2)--cycle;
\draw (L1)--(B1)--(A2)--(C2);
\draw(B3)--(B4);

\fill [red] (A0) circle (3pt);
\fill [red] (B4) circle (3pt);

\fill [black] (A2) circle (2pt);

\fill [black] (B1) circle (2pt);
\fill [black] (B3) circle (2pt);

\fill [black] (C2) circle (2pt);

\fill [black] (L1) circle (2pt);

\end{tikzpicture} &

\begin{tikzpicture}[scale=0.5]

\coordinate (A0) at (0,0);
\coordinate (A2) at ({sqrt(3)},0);

\coordinate (B1) at (0.5*{sqrt(3)}, 0.5);
\coordinate (B3) at (1.5*{sqrt(3)}, 0.5);
\coordinate (B4) at (1.5*{sqrt(3)}+1, 0.5);
\coordinate (B5) at (1.5*{sqrt(3)}+2, 0.5);

\coordinate (C2) at ({sqrt(3)},1);

\coordinate (L1) at (0.5*{sqrt(3)}, -0.5);

\filldraw[blue!20] (A0)--(L1)--(B3)--(C2)--cycle;

\draw (A0)--(L1)--(B3)--(C2)--cycle;
\draw (L1)--(B1)--(A2)--(C2);
\draw(B3)--(B5);

\fill [red] (A0) circle (3pt);
\fill [red] (B5) circle (3pt);

\fill [black] (A2) circle (2pt);

\fill [black] (B1) circle (2pt);
\fill [black] (B3) circle (2pt);
\fill [black] (B4) circle (2pt);

\fill [black] (C2) circle (2pt);

\fill [black] (L1) circle (2pt);

\end{tikzpicture} &

\begin{tikzpicture}[scale=0.5]

\coordinate (A0) at (0,0);
\coordinate (A2) at ({sqrt(3)},0);

\coordinate (B1) at (0.5*{sqrt(3)}, 0.5);
\coordinate (B3) at (1.5*{sqrt(3)}, 0.5);
\coordinate (B4) at (1.5*{sqrt(3)}+1, 0.5);
\coordinate (B5) at (1.5*{sqrt(3)}+2, 0.5);
\coordinate (B6) at (1.5*{sqrt(3)}+3, 0.5);

\coordinate (C2) at ({sqrt(3)},1);

\coordinate (L1) at (0.5*{sqrt(3)}, -0.5);

\filldraw[blue!20] (A0)--(L1)--(B3)--(C2)--cycle;

\draw (A0)--(L1)--(B3)--(C2)--cycle;
\draw (L1)--(B1)--(A2)--(C2);
\draw(B3)--(B6);

\fill [red] (A0) circle (3pt);
\fill [red] (B6) circle (3pt);

\fill [black] (A2) circle (2pt);

\fill [black] (B1) circle (2pt);
\fill [black] (B3) circle (2pt);
\fill [black] (B4) circle (2pt);
\fill [black] (B5) circle (2pt);

\fill [black] (C2) circle (2pt);

\fill [black] (L1) circle (2pt);

\end{tikzpicture} &

\begin{tikzpicture}[scale=0.5]

\coordinate (A0) at (0,0);
\coordinate (A2) at ({sqrt(3)},0);
\coordinate (A4) at (2*{sqrt(3)},0);

\coordinate (B1) at (0.5*{sqrt(3)}, 0.5);
\coordinate (B3) at (1.5*{sqrt(3)}, 0.5);
\coordinate (B5) at (2.5*{sqrt(3)}, 0.5);

\coordinate (C2) at ({sqrt(3)},1);
\coordinate (C4) at (2*{sqrt(3)},1);

\coordinate (L1) at (0.5*{sqrt(3)}, -0.5);

\filldraw[blue!20] (A0)--(L1)--(B3)--(A4)--(B5)--(C4)--(B3)--(C2)--cycle;

\draw (A0)--(L1)--(B3)--(A4)--(B5)--(C4)--(B3)--(C2)--cycle;
\draw (L1)--(B1)--(A2)--(C2);
\draw (A4)--(C4);

\fill [red] (A0) circle (3pt);
\fill [red] (B5) circle (3pt);

\fill [black] (A2) circle (2pt);
\fill [black] (A4) circle (2pt);

\fill [black] (B1) circle (2pt);
\fill [black] (B3) circle (2pt);

\fill [black] (C2) circle (2pt);
\fill [black] (C4) circle (2pt);

\fill [black] (L1) circle (2pt);

\end{tikzpicture} \\

$N16$ & $N17$ & $N18$ & $N19$ & $N20$\\

\begin{tikzpicture}[scale=0.5]

\coordinate (A0) at (0,0);
\coordinate (A2) at ({sqrt(3)},0);
\coordinate (A4) at (2*{sqrt(3)},0);

\coordinate (B1) at (0.5*{sqrt(3)}, 0.5);
\coordinate (B3) at (1.5*{sqrt(3)}, 0.5);
\coordinate (B5) at (2.5*{sqrt(3)}, 0.5);
\coordinate (B6) at (2.5*{sqrt(3)}+1, 0.5);

\coordinate (C2) at ({sqrt(3)},1);
\coordinate (C4) at (2*{sqrt(3)},1);

\coordinate (L1) at (0.5*{sqrt(3)}, -0.5);

\filldraw[blue!20] (A0)--(L1)--(B3)--(A4)--(B5)--(C4)--(B3)--(C2)--cycle;

\draw (A0)--(L1)--(B3)--(A4)--(B5)--(C4)--(B3)--(C2)--cycle;
\draw (L1)--(B1)--(A2)--(C2);
\draw (A4)--(C4);
\draw (B5)--(B6);

\fill [red] (A0) circle (3pt);
\fill [red] (B6) circle (3pt);

\fill [black] (A2) circle (2pt);
\fill [black] (A4) circle (2pt);

\fill [black] (B1) circle (2pt);
\fill [black] (B3) circle (2pt);
\fill [black] (B5) circle (2pt);

\fill [black] (C2) circle (2pt);
\fill [black] (C4) circle (2pt);

\fill [black] (L1) circle (2pt);

\end{tikzpicture} &

\begin{tikzpicture}[scale=0.5]

\coordinate (A0) at (0,0);
\coordinate (A2) at ({sqrt(3)},0);
\coordinate (A4) at (2*{sqrt(3)},0);

\coordinate (B1) at (0.5*{sqrt(3)}, 0.5);
\coordinate (B3) at (1.5*{sqrt(3)}, 0.5);
\coordinate (B5) at (2.5*{sqrt(3)}, 0.5);

\coordinate (C2) at ({sqrt(3)},1);
\coordinate (C4) at (2*{sqrt(3)},1);
\coordinate (C6) at (3*{sqrt(3)},1);

\coordinate (D5) at (2.5*{sqrt(3)}, 1.5);

\coordinate (L1) at (0.5*{sqrt(3)}, -0.5);

\filldraw[blue!20] (A0)--(L1)--(B3)--(A4)--(C6)--(D5)--(B3)--(C2)--cycle;

\draw (A0)--(L1)--(B3)--(A4)--(C6)--(D5)--(B3)--(C2)--cycle;
\draw (L1)--(B1)--(A2)--(C2);
\draw (A4)--(C4)--(B5)--(D5);

\fill [red] (A0) circle (3pt);
\fill [red] (C6) circle (3pt);

\fill [black] (A2) circle (2pt);
\fill [black] (A4) circle (2pt);

\fill [black] (B1) circle (2pt);
\fill [black] (B3) circle (2pt);
\fill [black] (B5) circle (2pt);

\fill [black] (C2) circle (2pt);
\fill [black] (C4) circle (2pt);

\fill [black] (D5) circle (2pt);

\fill [black] (L1) circle (2pt);

\end{tikzpicture} &

\begin{tikzpicture}[scale=0.5]

\coordinate (A0) at (0,0);
\coordinate (A2) at ({sqrt(3)},0);

\coordinate (B1) at (0.5*{sqrt(3)}, 0.5);
\coordinate (B3) at (1.5*{sqrt(3)}, 0.5);

\coordinate (C2) at ({sqrt(3)},1);
\coordinate (C4) at (2*{sqrt(3)},1);

\coordinate (D3) at (1.5*{sqrt(3)}, 1.5);

\coordinate (L1) at (0.5*{sqrt(3)}, -0.5);

\filldraw[blue!20] (A0)--(L1)--(C4)--(D3)--cycle;

\draw (A0)--(L1)--(C4)--(D3)--cycle;
\draw (L1)--(B1)--(A2)--(C2)--(B3)--(D3);

\fill [red] (A0) circle (3pt);
\fill [red] (C4) circle (3pt);

\fill [black] (A2) circle (2pt);

\fill [black] (B1) circle (2pt);
\fill [black] (B3) circle (2pt);

\fill [black] (C2) circle (2pt);

\fill [black] (D3) circle (2pt);

\fill [black] (L1) circle (2pt);

\end{tikzpicture} &

\begin{tikzpicture}[scale=0.5]

\coordinate (A0) at (0,0);
\coordinate (A2) at ({sqrt(3)},0);

\coordinate (B1) at (0.5*{sqrt(3)}, 0.5);
\coordinate (B3) at (1.5*{sqrt(3)}, 0.5);

\coordinate (C2) at ({sqrt(3)},1);
\coordinate (C4) at (2*{sqrt(3)},1);
\coordinate (C5) at (2*{sqrt(3)}+1,1);

\coordinate (D3) at (1.5*{sqrt(3)}, 1.5);

\coordinate (L1) at (0.5*{sqrt(3)}, -0.5);

\filldraw[blue!20] (A0)--(L1)--(C4)--(D3)--cycle;

\draw (A0)--(L1)--(C4)--(D3)--cycle;
\draw (L1)--(B1)--(A2)--(C2)--(B3)--(D3);
\draw (C4)--(C5);

\fill [red] (A0) circle (3pt);
\fill [red] (C5) circle (3pt);

\fill [black] (A2) circle (2pt);

\fill [black] (B1) circle (2pt);
\fill [black] (B3) circle (2pt);

\fill [black] (C2) circle (2pt);
\fill [black] (C4) circle (2pt);

\fill [black] (D3) circle (2pt);

\fill [black] (L1) circle (2pt);

\end{tikzpicture} &

\begin{tikzpicture}[scale=0.5]

\coordinate (A0) at (0,0);
\coordinate (A2) at ({sqrt(3)},0);

\coordinate (B1) at (0.5*{sqrt(3)}, 0.5);
\coordinate (B3) at (1.5*{sqrt(3)}, 0.5);

\coordinate (C2) at ({sqrt(3)},1);
\coordinate (C4) at (2*{sqrt(3)},1);
\coordinate (C5) at (2*{sqrt(3)}+1,1);
\coordinate (C6) at (2*{sqrt(3)}+2,1);

\coordinate (D3) at (1.5*{sqrt(3)}, 1.5);

\coordinate (L1) at (0.5*{sqrt(3)}, -0.5);

\filldraw[blue!20] (A0)--(L1)--(C4)--(D3)--cycle;

\draw (A0)--(L1)--(C4)--(D3)--cycle;
\draw (L1)--(B1)--(A2)--(C2)--(B3)--(D3);
\draw (C4)--(C6);

\fill [red] (A0) circle (3pt);
\fill [red] (C6) circle (3pt);

\fill [black] (A2) circle (2pt);

\fill [black] (B1) circle (2pt);
\fill [black] (B3) circle (2pt);

\fill [black] (C2) circle (2pt);
\fill [black] (C4) circle (2pt);
\fill [black] (C5) circle (2pt);

\fill [black] (D3) circle (2pt);

\fill [black] (L1) circle (2pt);

\end{tikzpicture} \\

$N21$ & $N22$ & $N23$ & $N24$ & $N25$\\

\begin{tikzpicture}[scale=0.5]

\coordinate (A0) at (0,0);
\coordinate (A2) at ({sqrt(3)},0);

\coordinate (B1) at (0.5*{sqrt(3)}, 0.5);
\coordinate (B3) at (1.5*{sqrt(3)}, 0.5);
\coordinate (B5) at (2.5*{sqrt(3)}, 0.5);

\coordinate (C2) at ({sqrt(3)},1);
\coordinate (C4) at (2*{sqrt(3)},1);
\coordinate (C6) at (3*{sqrt(3)},1);

\coordinate (D3) at (1.5*{sqrt(3)}, 1.5);
\coordinate (D5) at (2.5*{sqrt(3)}, 1.5);

\coordinate (L1) at (0.5*{sqrt(3)}, -0.5);

\filldraw[blue!20] (A0)--(L1)--(C4)--(B5)--(C6)--(D5)--(C4)--(D3)--cycle;

\draw (A0)--(L1)--(C4)--(B5)--(C6)--(D5)--(C4)--(D3)--cycle;
\draw (L1)--(B1)--(A2)--(C2)--(B3)--(D3);
\draw (B5)--(D5);

\fill [red] (A0) circle (3pt);
\fill [red] (C6) circle (3pt);

\fill [black] (A2) circle (2pt);

\fill [black] (B1) circle (2pt);
\fill [black] (B3) circle (2pt);
\fill [black] (B5) circle (2pt);

\fill [black] (C2) circle (2pt);
\fill [black] (C4) circle (2pt);

\fill [black] (D3) circle (2pt);
\fill [black] (D5) circle (2pt);

\fill [black] (L1) circle (2pt);

\end{tikzpicture} &

\begin{tikzpicture}[scale=0.5]

\coordinate (A0) at (0,0);
\coordinate (A2) at ({sqrt(3)},0);
\coordinate (A4) at (2*{sqrt(3)},0);

\coordinate (B1) at (0.5*{sqrt(3)},0.5);
\coordinate (B3) at (1.5*{sqrt(3)},0.5);
\coordinate (B5) at (2.5*{sqrt(3)},0.5);

\coordinate (C2) at ({sqrt(3)},1);
\coordinate (C4) at (2*{sqrt(3)},1);

\coordinate (L1) at (0.5*{sqrt(3)},-0.5);
\coordinate (L3) at (1.5*{sqrt(3)},-0.5);

\filldraw[blue!20] (A0)--(L1)--(A2)--(L3)--(B5)--(C4)--(B3)--(C2)--cycle;

\draw (A0)--(L1)--(A2)--(L3)--(B5)--(C4)--(B3)--(C2)--cycle;
\draw (L1)--(B1)--(A2)--(C2);
\draw (A2)--(B3);
\draw (L3)--(B3)--(A4)--(C4);

\fill [red] (A0) circle (3pt);
\fill [red] (B5) circle (3pt);

\fill [black] (A2) circle (2pt);
\fill [black] (A4) circle (2pt);

\fill [black] (B1) circle (2pt);
\fill [black] (B3) circle (2pt);

\fill [black] (C2) circle (2pt);
\fill [black] (C4) circle (2pt);

\fill [black] (L1) circle (2pt);
\fill [black] (L3) circle (2pt);

\end{tikzpicture} &

\begin{tikzpicture}[scale=0.5]

\coordinate (A0) at (0,0);
\coordinate (A2) at ({sqrt(3)},0);
\coordinate (A4) at (2*{sqrt(3)},0);

\coordinate (B1) at (0.5*{sqrt(3)},0.5);
\coordinate (B3) at (1.5*{sqrt(3)},0.5);
\coordinate (B5) at (2.5*{sqrt(3)},0.5);
\coordinate (B6) at (2.5*{sqrt(3)}+1,0.5);

\coordinate (C2) at ({sqrt(3)},1);
\coordinate (C4) at (2*{sqrt(3)},1);

\coordinate (L1) at (0.5*{sqrt(3)},-0.5);
\coordinate (L3) at (1.5*{sqrt(3)},-0.5);

\filldraw[blue!20] (A0)--(L1)--(A2)--(L3)--(B5)--(C4)--(B3)--(C2)--cycle;

\draw (A0)--(L1)--(A2)--(L3)--(B5)--(C4)--(B3)--(C2)--cycle;
\draw (L1)--(B1)--(A2)--(C2);
\draw (A2)--(B3);
\draw (L3)--(B3)--(A4)--(C4);
\draw (B5)--(B6);

\fill [red] (A0) circle (3pt);
\fill [red] (B6) circle (3pt);

\fill [black] (A2) circle (2pt);
\fill [black] (A4) circle (2pt);

\fill [black] (B1) circle (2pt);
\fill [black] (B3) circle (2pt);
\fill [black] (B5) circle (2pt);

\fill [black] (C2) circle (2pt);
\fill [black] (C4) circle (2pt);

\fill [black] (L1) circle (2pt);
\fill [black] (L3) circle (2pt);

\end{tikzpicture} &

\begin{tikzpicture}[scale=0.5]

\coordinate (A0) at (0,0);
\coordinate (A2) at ({sqrt(3)},0);
\coordinate (A4) at (2*{sqrt(3)},0);

\coordinate (B1) at (0.5*{sqrt(3)},0.5);
\coordinate (B3) at (1.5*{sqrt(3)},0.5);
\coordinate (B5) at (2.5*{sqrt(3)},0.5);

\coordinate (C2) at ({sqrt(3)},1);
\coordinate (C4) at (2*{sqrt(3)},1);
\coordinate (C6) at (3*{sqrt(3)},1);

\coordinate (D5) at (2.5*{sqrt(3)},1.5);

\coordinate (L1) at (0.5*{sqrt(3)},-0.5);
\coordinate (L3) at (1.5*{sqrt(3)},-0.5);

\filldraw[blue!20] (A0)--(L1)--(A2)--(L3)--(C6)--(D5)--(B3)--(C2)--cycle;

\draw (A0)--(L1)--(A2)--(L3)--(C6)--(D5)--(B3)--(C2)--cycle;
\draw (L1)--(B1)--(A2)--(C2);
\draw (A2)--(B3);
\draw (L3)--(B3)--(A4)--(C4)--(B5)--(D5);

\fill [red] (A0) circle (3pt);
\fill [red] (C6) circle (3pt);

\fill [black] (A2) circle (2pt);
\fill [black] (A4) circle (2pt);

\fill [black] (B1) circle (2pt);
\fill [black] (B3) circle (2pt);
\fill [black] (B5) circle (2pt);

\fill [black] (C2) circle (2pt);
\fill [black] (C4) circle (2pt);

\fill [black] (D5) circle (2pt);

\fill [black] (L1) circle (2pt);
\fill [black] (L3) circle (2pt);

\end{tikzpicture} \\
$N26$ & $N27$ & $N28$ & $N29$ & \\
\end{tabular}

\end{center}

\caption{The necklace galleries in \cite{ElMc02}. Necklace galleries are  built by identifying the large red outer vertices.\label{fig:necklace}}

\end{figure}
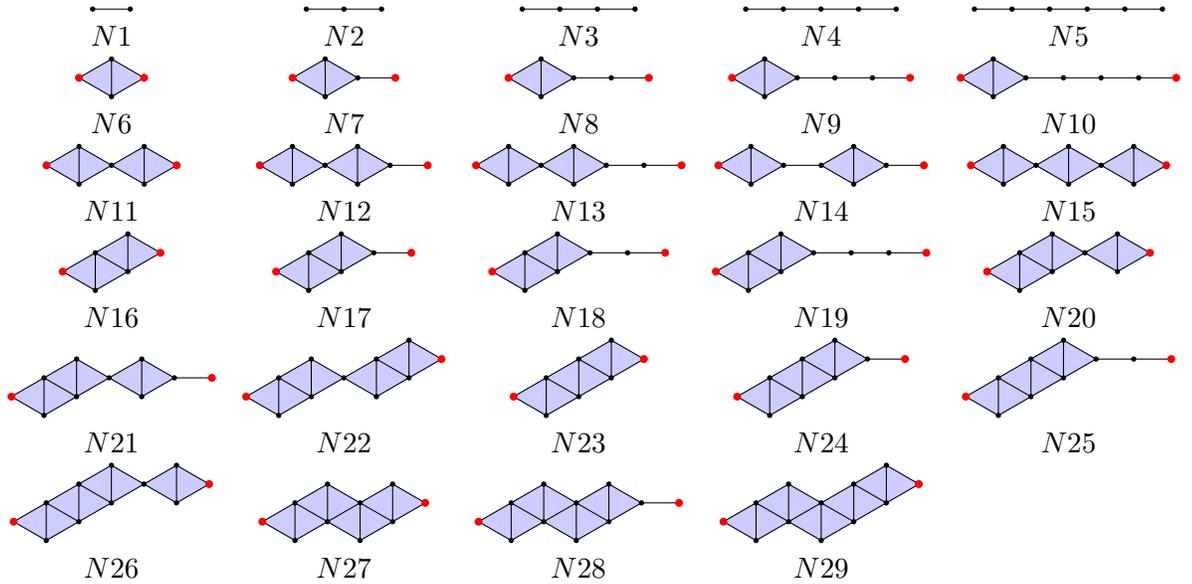

\noindent In  \cite{ElMc02}, Elder and McCammond  go on to give a specific list of configurations that must be avoided in links of a vertices of a non positively curved simplicial $3$-complex.  They broke these configuration into three types; annular galleries (Figure \ref{fig:annular}), m\"{o}bius galleries (Figure \ref{fig:mobius}), and necklace galleries (Figure \ref{fig:necklace}).

\begin{main}[Full $\Rightarrow$ NPC] Let $L$ be a full subcomplex of a $\cat(0)$ simplicial complex $K$. Then $L$ is nonpositively curved.
\end{main}

\begin{proof}  It suffices to prove that full subcomplexes avoid all of the annular, m\"{o}bius, and necklace galleries determined in \cite{ElMc02}. Each configuration will be referred to by its figure label. If a configuration occurs in $L$, it must be trivial in $K$. In other words, a configuration $c$ contained $lk_{L}(v)$ must be spanned by a spherical disk in $lk_K(v)$. 

Suppose $\alpha \subseteq lk_L(v)$ is a combinatorial triangle, pentagon, or square. By Theorem \ref{prop:linksnotempty}, $\alpha$ is filled by a disk $D$ in $lk_K(v)$ with no internal vertices. Since $lk_L(v)$ is full in $lk_K(v)$, $D \subseteq lk_L(v)$. Thus $lk_L(v)$ avoids configurations containing combinatorial paths of length less than six. This eliminates the annular galleries except $A8$ and $A9$, all of the m\"{o}bius galleries, and the necklace galleries except $N10$, $N13$, $N14$, $N15$, $N19$, $N21$, $N22$, $N25$, $N26$, $N28$ and $N29$. Each of these galleries contain $N10$ as as subdiagram (see Figure \ref{fig:n10inn25}). Thus eliminating $N10$ eliminates each of the remaining configurations.

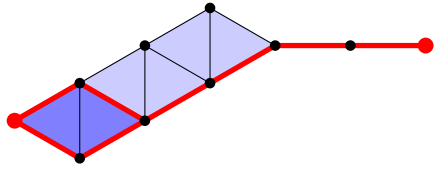
\begin{figure}

\begin{tikzpicture}

\coordinate (A0) at (0,0);
\coordinate (A2) at ({sqrt(3)},0);

\coordinate (B1) at (0.5*{sqrt(3)}, 0.5);
\coordinate (B3) at (1.5*{sqrt(3)}, 0.5);

\coordinate (C2) at ({sqrt(3)},1);
\coordinate (C4) at (2*{sqrt(3)},1);
\coordinate (C5) at (2*{sqrt(3)}+1,1);
\coordinate (C6) at (2*{sqrt(3)}+2,1);

\coordinate (D3) at (1.5*{sqrt(3)}, 1.5);

\coordinate (L1) at (0.5*{sqrt(3)}, -0.5);

\filldraw[blue!20] (A0)--(L1)--(C4)--(D3)--cycle;

\filldraw[blue!50] (A0)--(L1)--(A2)--(B1)--cycle;

\draw (A0)--(L1)--(C4)--(D3)--cycle;
\draw (L1)--(B1)--(A2)--(C2)--(B3)--(D3);
\draw (C4)--(C6);

\draw[red, line width=2pt] (A0)--(L1)--(A2)--(B1)--cycle;

\draw[red, line width=2pt] (A2)--(C4)--(C6);

\fill [red] (A0) circle (3pt);
\fill [red] (C6) circle (3pt);

\fill [black] (A2) circle (2pt);

\fill [black] (B1) circle (2pt);
\fill [black] (B3) circle (2pt);

\fill [black] (C2) circle (2pt);
\fill [black] (C4) circle (2pt);
\fill [black] (C5) circle (2pt);

\fill [black] (D3) circle (2pt);

\fill [black] (L1) circle (2pt);

\end{tikzpicture}

\caption{\label{fig:n10inn25} A copy of $N10$ in $N25$.}

\end{figure}

\begin{figure}

\begin{tikzpicture}

\coordinate (A0) at (0,0);
\coordinate (A2) at ({sqrt(3)},0);
\coordinate (A3) at ({sqrt(3)}+1,0);
\coordinate (A4) at ({sqrt(3)}+2,0);
\coordinate (A5) at ({sqrt(3)}+3,0);
\coordinate (A6) at ({sqrt(3)}+4,0);

\coordinate (B1) at (0.5*{sqrt(3)}, 0.5);

\coordinate (L1) at (0.5*{sqrt(3)}, -0.5);

\filldraw[blue!20] (A0)--(L1)--(A2)--(B1)--cycle;

\draw (A0)--(L1)--(A2)--(B1)--cycle;
\draw (A2)--(A6);

\draw (L1)--(B1);

\fill [black] (A0) circle (2pt) node [anchor=east] {$v_0$};
\fill [black] (A2) circle (2pt) node [anchor=south] {$v_2$};
\fill [black] (A3) circle (2pt) node [anchor=south] {$v_3$};
\fill [black] (A4) circle (2pt) node [anchor=south] {$v_4$};
\fill [black] (A5) circle (2pt) node [anchor=south] {$v_5$};
\fill [black] (A6) circle (2pt) node [anchor=south] {$v_6 = v_0$};

\fill [black] (B1) circle (2pt) node [anchor= south] {$v_1$};
\fill [black] (L1) circle (2pt);

\fill [red] (A0) circle (3pt);
\fill [red] (A6) circle (3pt);

\end{tikzpicture}

\caption{\label{fig:n10label}The vertex labeling of the upper path in $N10$ used in the proof of the main theorem.}

\end{figure}
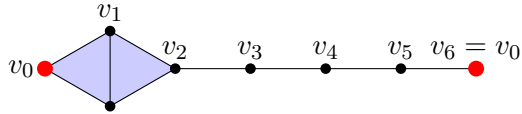

If $N10$ is in $lk_L(v)$ for a vertex $v$ in $L$, it must be trivial in $lk_K(v)$. There are two distinct non trivial combinatorial paths in $N10$; one along the lower edges of the triangle pair and one along the upper edges. The two cases are symmetric, so without loss of generality assume the path along the upper edge of the move is filled by a minimal spanning disk $D$ in $lk_K(v)$. Label in the vertices on $N10$ $v_0$ through $v_6=v_0$ (see Figure \ref{fig:n10label}). The vertex $v_1$ must lie on at least four triangles in $D$.  This implies that $deg(v_1) \geq 5$. Since $lk_L(v)$ is full in $lk_K(v)$, there are no edges in $D$ connecting two non consecutive vertices of the gallery. This implies $deg(v_i) \geq 3$ for $i \neq 1$. Thus \begin{equation*}\sum_{i=1}^6 (4-deg(v_i))= (4-deg(v_1)) + \sum_{i=2}^6(4- deg(v_i))\\ \leq -1 + 5(1) = 4,\end{equation*}  a contradiction to Combinatorial Gau\ss-Bonnet.

\end{proof}

\begin{cor}[Full in NPC $\Rightarrow$ NPC]\label{cor:fullnpc} A full subcomplex of a nonpositively curved simplicial $3$-complex is nonpostively curved.
\end{cor}

\begin{proof} Let $L$ be a full subcomplex of a nonpositively curved simplicial $3$-complex $K$. Then the universal cover $\widetilde{L}$ is a full subcomplex of the universal cover $\widetilde{K}$. $\widetilde{K}$ is $\cat(0)$,  so $\widetilde{L}$ is nonpositively curved by our main theorem. As curvature is a local condition, this implies $L$ is also nonpositively curved.
\end{proof}




\section{Higher Dimensions}\label{sec:highdim}


Unfortunately the main theorem does not generalized to higher dimensional simplicial complexes. While higher dimensional complexes satisfy the no empty triangle and square conditions, they might have empty pentagons.

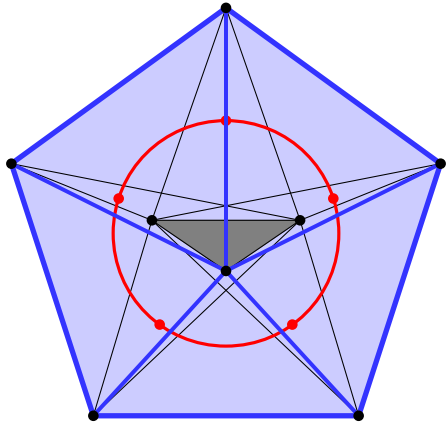
\begin{figure} \label{fig:narwahl0}

\begin{tikzpicture}

\coordinate (O) at (0,0);
\coordinate (P0) at (90: 3);
\coordinate (P1) at (90+72: 3);
\coordinate (P2) at (90+2*72: 3);
\coordinate (P3) at (90+3*72: 3);
\coordinate (P4) at (90+4*72: 3);

\coordinate (T0) at (270 + 100: 1);
\coordinate (T1) at (270 - 100: 1);
\coordinate (T2) at (270: 0.5);

\coordinate (L0) at (90: 1.5);
\coordinate (L1) at (90+72: 1.5);
\coordinate (L2) at (90+2*72: 1.5);
\coordinate (L3) at (90+3*72: 1.5);
\coordinate (L4) at (90+4*72: 1.5);

\filldraw[blue!20] (P0)--(P1)--(P2)--(P3)--(P4)--cycle;

\draw[blue!80, line width=2 pt] (P0)--(P1)--(P2)--(P3)--(P4)--cycle;

\filldraw[black!50] (T0)--(T1)--(T2)--cycle;

\draw (T0)--(T1)--(T2)--cycle;

\draw (T0)--(P0);
\draw (T0)--(P1);
\draw (T0)--(P2);
\draw (T0)--(P3);
\draw (T0)--(P4);

\draw (T1)--(P0);
\draw (T1)--(P1);
\draw (T1)--(P2);
\draw (T1)--(P3);
\draw (T1)--(P4);

\draw[red, very thick] (O) circle (1.5);

\fill[red] (L0) circle (2pt);
\fill[red] (L1) circle (2pt);
\fill[red] (L2) circle (2pt);
\fill[red] (L3) circle (2pt);
\fill[red] (L4) circle (2pt);

\draw[blue!80, line width=1.5 pt] (T2)--(P0);
\draw[blue!80, line width=1.5 pt] (T2)--(P1);
\draw[blue!80, line width=1.5 pt] (T2)--(P2);
\draw[blue!80, line width=1.5 pt] (T2)--(P3);
\draw[blue!80, line width=1.5 pt] (T2)--(P4);


\fill[black] (P0) circle (2pt);
\fill[black] (P1) circle (2pt);
\fill[black] (P2) circle (2pt);
\fill[black] (P3) circle (2pt);
\fill[black] (P4) circle (2pt);

\fill[black] (T0) circle (2pt);
\fill[black] (T1) circle (2pt);
\fill[black] (T2) circle (2pt);

\end{tikzpicture}

\caption{The complex described in Example \ref{exmp:baddisk} for $n=4$. \label{fig:baddisk}}

\end{figure}

\begin{exmp}[High Dimensional Empty Pentagon]\label{exmp:baddisk} Let $\sigma$ be an $(n-2)$-dimensional simplex for $n \geq 4$, and take $K$ to be the join of  $\sigma$ together with closed cycle of five distinct vertices $\alpha = [v_0,v_2,...,v_5=v_0]$. Then $K$ is the union of five distinct $n$-simplices arranged cyclically around $\sigma$ (see Figure \ref{fig:baddisk}).  The curvature of $K$ depends only on $lk_K(\sigma)$. The link is a regular metric graph consisting of a single cycle $c$ with five edges. Each edge has length $\arccos(\frac{1}{n})$, giving $\ell(c) = 5\arccos(\frac{1}{n}) > 2\pi$ for $n \geq 4$.  $K$ is also simply connected, so $K$ is $\cat(0)$.  There are $n-1$ minimal disks spanning $\alpha$ in $K$, one through each of the vertices of $\sigma$. Each disk is an example of a full subcomplex that is positively curved; the disks consist of five triangles around a interior vertex. \end{exmp}


\nocite{BrHa99}

\nocite{Bo90}


\nocite{Gr87}


Excluding configurations like those described in Example \ref{exmp:baddisk} give complexes which fall under Januszkiewicz and \'{S}wi\c{a}tkowski's theory of simplicial non positive curvature \cite{JaSw06}. The simplicial nonpositive curvature condition is a combinatorial condition on simplicial complexes that is similar to the $\cat(0)$ condition.

\begin{defn}[Simplicial Nonpositive Curvature] The \emph{combinatorial link} of a simplex $\sigma$, denoted $clk_{K}(\sigma)$ is the union of all simplicies $\tau$ of $K$ such that $\sigma \ast \tau$ is a simplex of $K$. A simplicial complex $K$ is \emph{$k$-large} if $K$ is flag and contains no empty $n$-gons for $n<k$. A simplicial complex $K$ satisfies the \emph{simplicial nonpositive curvature} condition, or \emph{$\textsc{SNPC}$},  if $clk_{K}(\sigma)$ is $6$-large  for each simplex $\sigma \subseteq K$. \end{defn}

The fact that full subcomplexes of $\textsc{SNPC}$ complexes inherit the $\textsc{SNPC}$ condition immediately follows from the definion \cite{JaSw06}.

\bibliography{researchbib}

\begin{thebibliography}{1}

\bibitem{Bo90}
B.~H. Bowditch.
\newblock Notes on {G}romov's hyperbolicity criterion for path-metric spaces.
\newblock In {\em Group theory from a geometrical viewpoint ({T}rieste, 1990)},
  pages 64--167. World Sci. Publ., River Edge, NJ, 1991.

\bibitem{Bo92}
B.~H. Bowditch.
\newblock Notes on locally {${\rm CAT}(1)$} spaces.
\newblock In {\em Geometric group theory ({C}olumbus, {OH}, 1992)}, volume~3 of
  {\em Ohio State Univ. Math. Res. Inst. Publ.}, pages 1--48. de Gruyter,
  Berlin, 1995.

\bibitem{BrHa99}
Martin~R. Bridson and Andre Haefliger.
\newblock {\em Metric spaces of non-positive curvature}.
\newblock Springer-Verlag Berlin Heidelberg, New York, NY, 1999.

\bibitem{Cr08}
Katherine Crowley.
\newblock Simplicial collapsibility, discrete {M}orse theory, and the geometry
  of nonpositively curved simplicial complexes.
\newblock {\em Geom. Dedicata}, 133:35--50, 2008.

\bibitem{ElMc02}
Murray Elder and Jon McCammond.
\newblock Curvature testing in 3-dimensional metric polyhedral complexes.
\newblock {\em Experiment. Math.}, 11(1):143--158, 2002.

\bibitem{ElMc04}
Murray Elder and Jon McCammond.
\newblock C{AT}(0) is an algorithmic property.
\newblock {\em Geom. Dedicata}, 107:25--46, 2004.

\bibitem{Gr87}
M.~Gromov.
\newblock Hyperbolic groups.
\newblock In {\em Essays in group theory}, volume~8 of {\em Math. Sci. Res.
  Inst. Publ.}, pages 75--263. Springer, New York, 1987.

\bibitem{JaSw06}
Tadeusz Januszkiewicz and Jacek Swiatkowski.
\newblock Simplicial nonpositive curvature.
\newblock {\em Publications Mathematiques de l'IHES}, 2006.
\newblock to appear.

\end{thebibliography}

\bibliographystyle{plain}

\end{document}